%% file: JLM_15.tex
\documentclass[12pt]{amsart}

\usepackage{latexsym}

\usepackage{amssymb}
\usepackage[centertags]{amsmath}
\usepackage{verbatim}
\usepackage{epic, eepic}

\newtheorem{theorem}{Theorem}
\newtheorem{lemma}{Lemma}
\newtheorem{claim}[equation]{Claim}
\newtheorem*{claim*}{Claim}
\newtheorem{corollary}[equation]{Corollary}

\newcommand{\J}{\mathcal J}

\renewcommand{\ge}{\geqslant}
\renewcommand{\le}{\leqslant}

\renewcommand{\phi}{\varphi}

\renewcommand{\#}{\operatorname{Card}}

\renewcommand{\epsilon}{\varepsilon}

\newcommand{\Ex}{\mathbb E}
\renewcommand{\Pr}[1]{\,\mathbb P\,\big\{\,#1\,\big\}\,}
\newcommand{\Pro}[1]{\,\mathbb P\,\Big\{\,#1\,\Big\}\,}

\renewcommand{\le}{\leqslant}
\renewcommand{\ge}{\geqslant}

\newcommand{\C}{\mathbb C}
\newcommand{\Z}{\mathbb Z}
\newcommand{\D}{\mathbb D}

\newcommand{\T}{\mathbb T}

\title[Large fluctuations of random
complex zeroes]{The Jancovici - Lebowitz - Manificat law \\ for
large fluctuations of random complex zeroes}

\author{F. Nazarov}
\address{F.N.: Mathematics Department \\
University of Wisconsin-Madison \\
480 Lincoln Dr., Madison WI 53706\\
USA}

\email{nazarov@math.wisc.edu}

\author{M. Sodin}
\address{M.S.: School of Mathematical Sciences\\
Tel Aviv University\\
Tel Aviv 69978\\
Israel}

\email{sodin@post.tau.ac.il}

\author{A. Volberg}
\address{A.V.: Department of Mathematics\\
Michigan State University\\
East Lansing, MI 48824\\
USA} \email{volberg@math.msu.edu}

\thanks{F.N. and A.V. are partially supported by
the National Science Foundation, DMS grant 0501067. M.S. is
partially supported by the Israel Science Foundation of the Israel
Academy of Sciences and Humanities, grants 357/04 and 171/07}

\date{}

\begin{document}

\begin{abstract}
Consider a Gaussian Entire Function
\[
f(z) = \sum_{k=0}^\infty \zeta_k \frac{z^k}{\sqrt{k!}} \, ,
\]
where $\zeta_0, \zeta_1, \dots $ are Gaussian i.i.d. complex
random variables. The zero set of this function is distribution
invariant with respect to the isometries of the complex plane. Let
$n(R)$ be the number of zeroes of $f$ in the disk of radius $R$. It
is easy to see that $\Ex n(R) = R^2$, and it is known that the
variance of $n(R)$ grows linearly with $R$ (Forrester and Honner).
We prove that, for every $ \alpha > 1/2$, the tail probability
$\Pr { | n(R)-R^2 | > R^\alpha }$ behaves as $\exp\left[ -R^{\phi
(\alpha)} \right]$ with some explicit piecewise linear function
$\phi (\alpha)$. For some special values of the parameter
$\alpha$, this law was found earlier by Sodin and Tsirelson, and
by Krishnapur.

In the context of charge fluctuations of a one-component Coulomb
system
of particles of one sign embedded into a uniform background of
another sign, a similar law was discovered some time ago by
Jancovici, Lebowitz
and Manificat.

\end{abstract}

\maketitle

\section{Introduction}\label{intro}
\numberwithin{equation}{section}

Consider the Fock-Bargmann space of the entire functions of one
complex variable that are square integrable with respect to the
measure $\tfrac1{\pi}e^{-|z|^2}\,dm(z)$, where $m$ is the Lebesgue
measure on $\C$. Let $f$ be a Gaussian function associated with
this space; i.e.,
\[
f(z) = \sum_{k\ge 0}\zeta_k e_k(z)
\]
where $\zeta_k$ are independent standard complex Gaussian random
variables (that is, the density of $\zeta_k$ on the complex plane
$\C$ is $\tfrac1\pi e^{-|w|^2}$), and $\{e_k\}$ is an orthonormal
basis in the Fock-Bargmann space. The Gaussian function $f$ does
not depend on the choice of the basis $\{e_k\}$, so usually one
takes the standard basis $\displaystyle e_k(z) =
\frac{z^k}{\sqrt{k!}}$, $k\in \Z_+$. In what follows, we call $f$
a {\em Gaussian Entire Function} (G.E.F., for short). G.E.F.
together with other similar models were introduced in the 90's in
the works of Bogomolny, Bohigas, Lebouef~\cite{BBL}, and
Hannay~\cite{H}.

A remarkable feature of the zero set $\mathcal Z_f = f^{-1}\{ 0
\}$ of a G.E.F. is its distribution invariance with respect to the
isometries of $\C$. The rotation invariance is obvious since the
distribution of the function $f$ is rotation invariant. The
translation invariance follows, for instance, from the fact that
the operators $(T_w g)(z) = g(w+z)e^{-z\overline{w}}
e^{-|w|^2/2}$, $w\in\C$, are unitary operators in the
Fock-Bargmann space, and therefore, if $f$ is a G.E.F., then $T_w
f$ is a G.E.F. as well (see section~2.2 below). It is worth
mentioning that by Calabi's rigidity \cite[Section~3]{S}, $f(z)$
together with its scalings $f(tz)$, $t>0$, are  the only Gaussian
entire functions with the distribution of zeroes invariant with
respect to the isometries of $\C$. See \cite[Part~I]{ST} for
further discussion.

Let $n(R) = \# \big\{ \mathcal Z_f \cap R\D\big\} $ be the number
of zeroes of $f$ in the disk of radius $R$. It is not hard to
check that the mean number of points of $\mathcal Z_f$ per unit
area equals $\tfrac1{\pi}$ (cf. Section~\ref{sect2.1.1}).
Therefore, $\Ex n(R) = R^2$. The asymptotics of the variance of
$n(R)$ was computed by Forrester and Honner in \cite{FH}: \[ \Ex
\big( n(R) - R^2 \big)^2 = cR + o(R), \qquad R\to\infty, \] with
an explicitly computed positive $c$. In \cite{SZ}, Shiffman and
Zelditch gave a different computation of the asymptotics of the
variance valid in a more general context. The normalized random
variables $\displaystyle \frac{n(R)-R^2}{\sqrt{\operatorname{Var}
n(R)}} $ converge in distribution to the standard Gaussian random
variable. This can be proven, for instance, by a suitable
modification of the argument used in~\cite[Part I]{ST}. In this
work, we describe the probabilities of large fluctuations of the
random variable $n(R) - R^2$.
\begin{theorem}\label{thm_main}
For every $ \alpha>\tfrac12 $ and every $\epsilon>0$,
\begin{equation}\label{eq_JLM}
e^{-R^{\phi (\alpha) + \epsilon}} < \Pr { | n(R)-R^2 |
> R^\alpha } < e^{ -R^{\phi (\alpha) - \epsilon}}
\end{equation}
for all sufficiently large $R>R_0(\alpha, \epsilon)$, where
\[
\phi (\alpha) =
\begin{cases}
2\alpha-1, & \tfrac12 < \alpha \le 1; \\
3\alpha-2, & 1 \le \alpha \le 2; \\
2\alpha, & \alpha \ge 2\,.
\end{cases}
\]
\end{theorem}

\medskip

In a different context of charge fluctuations of a one-component
Coulomb system of particles of one sign embedded into a uniform
background of the opposite sign, a similar law was discovered by
Jancovici, Lebowitz and Manificat in their physical paper
\cite{JLM}. Let us mention that it is known since Ginibre's
classical paper \cite{G} that the class of point processes
considered by Jancovici, Lebowitz and Manificat contains as a
special case the $N\to\infty$ limit of the eigenvalue point
process of the ensemble of $N\times N$ random matrices with
independent standard complex Gaussian entries. The resemblance
between the zeroes of G.E.F. and the eigenvalues of Ginibre's
ensemble was discussed both in the physical and the mathematical
literature.

\medskip

Now, let us return to the zeroes of G.E.F.\,. In some cases, the
estimate \eqref{eq_JLM} is known. As we have already mentioned, it
is known for $\alpha = \tfrac12$ when it follows from the
asymptotics of the variance and the asymptotic normality. In the
case $\alpha = 2$ it follows from a result of Sodin and
Tsirelson~\cite[Part III]{ST}, which says that {\em for each $R\ge
1$,}
\[
e^{-CR^4} \le  \Pr { |n(R)-R^2| > R^2 } \le e^{-c R^4}
\]
with some positive numerical constants $c$ and $C$.
In~\cite{Krishnapur}, Krishnapur considered the case $\alpha > 2$
and proved that in that case
\[
\Pr { n(R) > R^\alpha } =   e^{- (\frac{\alpha}2 -1)(1+o(1))
R^{2\alpha} \log R  }\,, \qquad R\to\infty\,.
\]
In the same work, he also proved the lower bound in the case
$1<\alpha<2$:
\[
\Pr { | n(R) - R^2 | > R^\alpha } \ge e^{- C R^{3\alpha-2}}\,.
\]
Using a certain development of his method, we'll get the lower
bound
\[
\Pr { | n(R) - R^2 | > R^\alpha } \ge e^{- CR^{2\alpha-1}}\,,
\qquad \frac12<\alpha<1\,.
\]
Apparently, in the case $\tfrac12< \alpha <2$,
the technique used in \cite[Part III]{ST} and \cite{Krishnapur}
does not allow one to treat the upper bounds in the law
\eqref{eq_JLM}, which require new ideas.

\bigskip\par\noindent{\bf Outline of the proof.} Let us sketch the
main ideas we use in the proof of Theorem~\ref{thm_main}.

\medskip\par\noindent{\bf 1.}
We denote by $\Delta_I \arg f$ the increment of the argument of a
G.E.F. $f$ over an arc $I\subset R\T$ oriented counterclockwise,
and set
\[
\delta (f, I) = \Delta_I \arg f - \Ex \Delta_I \arg f = \Delta_I
\arg f - |I|\cdot R\,.
\]
Then by the argument principle,
\[
2\pi (n(R) - R^2) = \delta(f, R\T)\,.
\]
Note that the random variable $\delta(f, I)$ is set-additive and
split the circumference $R\T$ into $\displaystyle N = 2\pi
\frac{R}{r}$ disjoint arcs $I_j$ of length $r$. Thus we need to
estimate the probability of the event
\[
\Omega_\alpha (R)= \left\{\, \Bigl| \sum_{j=1}^N \delta(f, I_j)
\Bigr| > 2\pi R^\alpha\right\}.
\]

\medskip\par\noindent{\bf 2.} Let us fix an arc $I$ of length $r$
and look more closely at the tails of the random variable $\delta
(f, I)$.
It is not hard to check that $\delta(f, I) = \delta (T_w f, I-w)$
where $w$ is the center of the arc
$I$ and $T_w f(z) = f(w+z) e^{-z\overline{w}}e^{-|w|^2/2}$. A
classical complex analysis argument shows that for any analytic
function $g$ in the disk $2r\D$ and any ``good" arc $\gamma\subset
r\D$ of length at most $r$, one has
\[
\big | \Delta_\gamma \arg g \big| \le C\,
\log\frac{\max_{2r\D}|g|}{\max_{r\D} |g|}\,,
\]
see Lemma~\ref{lemma6.2}. Then, estimating the probability that,
for a G.E.F. $g=T_wf$, the doubling exponent $\displaystyle
\log\frac{\max_{2r\D}|g|}{\max_{r\D} |g|} $ is large, we come up
with the tail estimate
\[
\Pr { \big| \delta(f, I) \big| > Mr^2 } \le \exp\Bigl(-
\frac{CM^2}{\log M}\, r^4\Bigr)\,, \qquad M\gg 1\,.
\]

\medskip\par\noindent{\bf 3.} Now, let us come back to the sum
$\displaystyle \sum_{j=1}^N \delta(f, I_j)$. The random variables
$\delta(f, I_j)$ are not independent, however in
\cite[Theorem~3.2]{NSV} we've introduced an ``almost independence
device'' that allows us to think about these random variables as
of independent ones, provided that the arcs $I_j$ are
well-separated from each other. Here we'll need a certain
extension of that result (Lemma~\ref{thm3.1} below).

\medskip\par\noindent{\bf 4.} To see how the almost independence
and
the tail estimate work, first, consider the case $1<\alpha<2$. We
split the circumference $R\T$ into $N$ disjoint arcs $\{I_j\}$ of
length $r$. In view of the tail estimate in item 2, we need to
distribute the total deviation $R^\alpha$ between these arcs in
such a way that the ``deviation per arc'' $R^\alpha/N$ is bigger
than $r^2$. Since $N\simeq \tfrac{R}{r}$, this leads to the choice
of $r$ comparable to $R^{\alpha-1}$.

Then we consider the event that for a {\em fixed subset} $J\subset
\{1, 2, \, ...\,, N\}$ and for every $j\in J$, one has $|\delta(f, I_j)| \ge m_j r^2$, where $m_j$ are some big positive integer
powers of $2$ that satisfy
\begin{equation}\label{eq*}
\sum_{j\in J} m_j r^2 \gtrsim R^\alpha\,.
\end{equation}
Then we choose a well-separated sub-collection of arcs $J'\subset
J$ that falls under the assumptions of the almost independence
lemma~\ref{thm3.1}.
This step weakens condition~\eqref{eq*} to
\[
\sum_{j\in J'} m_j^{3/2} r^2 \gtrsim R^\alpha\,,
\]
which still suffices for our purposes. Then regarding the random
variables $\delta
(f, I_j)$, $j\in J'$, as independent ones and using the tail
estimate for these
variables, we see that the probability of this event does not
exceed
\begin{multline*}
\exp\Bigl( -c \sum_{j\in J'} \frac{m_j^2}{\log m_j} r^4 \Bigr) \le
\exp\Bigl( - r^2 \sum_{j\in J'} m_j^{3/2} r^2  \Bigr) \\ \le
\exp\bigl( - c r^2 R^\alpha \bigr) \le \exp\bigl( - c_1 R^{3\alpha
-2} \bigr)\,.
\end{multline*}

To get the upper bound for the probability of the event
$\Omega_\alpha$, we need to take into account the number of
possible choices of the subset $J$ and of the numbers $m_j$. This
factor does not exceed $ 2^N (\log R)^N < e^{CR\log\log R}$ which
is not big enough to destroy our estimate.

\medskip\par\noindent{\bf 5.} Now, let us turn to the upper bound
in the case
$\tfrac12 < \alpha < 1$. We choose the arcs $I_j$ of length
$1$. To separate them from each other, we choose from this
collection $R^{1-\epsilon}$ arcs $\{I_j\}_{j\in J}$ separated by
$R^\epsilon$ and such that
\[
\Bigl| \sum_{j\in J} \delta(f, I_j) \Bigr| > R^{\alpha -
\epsilon}\,.
\]
For these arcs, the random variables $\delta(f, I_j)$ behave like
independent ones, and since their tails have a fast decay, we can
apply to them the classical Bernstein inequality
(Lemma~\ref{lemma2.3}), which yields
\[
\Pro { \Big| \sum_{j\in J} \delta(f, I_j) \Big| > R^{\alpha -
\epsilon} }  \le C \exp\Bigl( - \frac{c(R^{\alpha -
\epsilon})^2}{\# J} \Bigr) = C \exp \bigl( -c R^{2\alpha -1 -
\epsilon} \bigr)\,.
\]

\smallskip\par\noindent{\bf 6.} To get the lower bound for the
probability
of $\Omega_\alpha$ in the case $\tfrac12<\alpha<1$, we introduce
an auxiliary Gaussian Taylor series
\[
g(z) = \sum_{k=0}^\infty \zeta_k a_k \frac{z^k}{\sqrt{k!}}
\]
where $\zeta_k$ are independent standard complex Gaussian random
variables, and
\[
a_k =
\begin{cases}
\sqrt{1-R^{\alpha-1}},& R^2 + R < k < R^2 + 2R\,;
\\
\sqrt{1+R^{\alpha-1}},& R^2-2R < k < R^2-R\,;
\\
1,& {\rm otherwise}\,.
\end{cases}
\]
It is not difficult to check that for some absolute $c>0$, the
probability that the function $g$ has at most $R^2-cR^\alpha$
zeroes in the disk $R\D$ is not exponentially small (more
precisely, it cannot be less than $c R^{-2+\alpha}$).

Now, let $\gamma $ be the standard Gaussian measure in the space
$\C^\infty$; i.e., the product of countably many copies of standard
complex Gaussian measures on $\C$, and let $\gamma_a$ be
another Gaussian measure on $\C^\infty$ which is the product of
complex Gaussian measures $\gamma_{a_k}$ on $\C$ with
variances $a_k^2$. Let $E\subset \C^\infty $ be the set of
coefficients $\eta_k$ such that the Taylor series $\sum_{k\ge
0}\eta_k\frac{z^k}{\sqrt{k!}}$ converges in $\C$ and has at most
$R^2-c R^{\alpha}$ zeroes in $R\D$. Then
\[
\gamma_a (E) \ge c R^{-2+\alpha}\,,
\]
while the quantity $\Pr { n(R)\le R^2 - c R^\alpha}$ we are
interested in equals $\gamma (E)$. Thus, it remains to compare
$\gamma (E)$ with $\gamma_a (E)$, and a more or less
straightforward computation finishes the job.

\bigskip

Let us mention that in the range $\tfrac12 < \alpha <1 $, the
exponent $\phi (\alpha) = 2\alpha-1$ is universal and seems to be
determined by the asymptotic normality at the endpoint
$\alpha=\tfrac12$. On the other hand, in the range $\alpha > 1$,
the law \eqref{eq_JLM} is not universal. To illustrate this point,
as in \cite[Parts I and II]{ST}, we consider random independent
perturbations of the lattice points. We fix the parameter $\nu
> 0$, and consider the random point set $ \{ \omega +
\zeta_{\omega}\}_{\omega\in\mathbb Z^2} $, where $\zeta_\omega$
are independent, identical, radially distributed random variables
with the tails  $\Pr { |\zeta_\omega|> t } $ decaying as $\exp (-t
^\nu)$ for $t\to\infty$. Set \[ n(R) = \# \{\omega\in \mathbb
Z^2\colon |\omega+\zeta_\omega|\le R \}.\] Then one can see that,
for every $ \alpha>\tfrac12 $ and every $\epsilon>0$,
\[
e^{-R^{\phi (\alpha, \nu) + \epsilon}} < \Pr { | n(R) - \pi R^2 |
> R^\alpha } < e^{-R^{\phi (\alpha, \nu) - \epsilon}}
\]
for all sufficiently large $R>R_0(\alpha, \epsilon)$ with
\[
\phi (\alpha, \nu) =
\begin{cases}
2\alpha-1, & \tfrac12 \le \alpha \le 1; \\
(\nu+1)\alpha - \nu, & 1 \le \alpha \le 2; \\
(\nu/2+1)\alpha, & \alpha \ge 2\,.
\end{cases}
\]
That is, the Jancovici-Lebowitz-Manificat law \eqref{eq_JLM}
corresponds to the case $\nu = 2$; i.e., to the lattice
perturbation with the Gaussian decay of the tails. It is
interesting to juxtapose this observation with the results from
\cite[Part II]{ST} about the matching between the zero set
$\mathcal Z_f$ and the lattice $\displaystyle \frac1{\sqrt{\pi}}
\Z^2$, and from \cite{NSV} about the gradient transportation of
the area measure to $\mathcal Z_f$.

\bigskip\par\noindent{\bf Convention about the constants.}
By $c$ and $C$ we denote positive numerical constants that appear
in the proofs. The constants denoted by $c$ are supposed to be
small
(in particular, they are always less than $1$), while the constants
denoted by $C$ are supposed to be big (they are always larger than
$1$). Within the proof of each lemma, we start a new
sequence of indices for these constants, and we never refer to
these constants after the corresponding proof is completed.

\medskip Notation $A \lesssim B$ and $ A \gtrsim B $ means that
there exist
positive numerical constants $C$ and $c$ such that $A\le C\cdot B$
and $ A \ge c \cdot B $ correspondingly. If $A\lesssim B$ and
$A\gtrsim B$
simultaneously, then we write $A\simeq B$. Notation $A \ll B$
stands for ``much less'' and means that $A \le c \cdot B$ with a
very small positive $c$; similarly, $A \gg B$ stands for ``much
larger'' and means that $A\ge C\cdot B$ with a very large positive
$C$.

\bigskip\par\noindent{\bf Acknowledgement.} We thank Manjunath
Krishnapur,
Yuval Peres, and Boris Tsirelson for very helpful discussions.

\numberwithin{equation}{subsection}

\section{Preliminaries}\label{sect2}

\subsection{A combinatorial lemma}

For $j, k
\in \{1, \, ...\,, N\}$, we set
\begin{multline*}
|j-k|_* = \min \left\{ |i-k|\colon i\equiv j \mod N \right\} \\
= \min\left\{ |j-k|, \ |j-k+N|, \ |j-k-N| \right\}\,.
\end{multline*}

\begin{lemma}\label{claim7.2}
Let $m_1$, ..., $m_N$ be non-negative integers. Then,
given $Q \ge 1$, there exists a subset $J'\subset \left\{ 1,\,
...\,, N\right\}$
such that
\[
|j-k|_* \ge Q (\sqrt{m_j} + \sqrt{m_k})\,, \qquad j, k \in J',
\quad j\ne k\,,
\]
and
\[
\sum_{j\in J} m_j \le 5Q \sum_{j\in J'} m_j^{3/2}\,.
\]
\end{lemma}

\par\noindent{\em Proof of Lemma~\ref{claim7.2}:} We build the set
$J'$ by an inductive construction. Choose $j_1\in
\{1, \, ...\,, N\}$ such that $ m_{j_1} = \max \big\{ m_j\colon
j\in \{1, \, ...\,, N\} \big\}$. Set
\[
J_1' = \{j_1\}, \quad J_1'' = \big\{j\colon 0<|j-j_1|_* <
2Q\sqrt{m_{j_1}} \big\}, \quad J_1 = J_1' \cup J_1'',
\]
and note that
\[
\sum_{j\in J_1} m_j \le \left( 4Q\sqrt{m_{j_1}} +1 \right) m_{j_1}
\le 5Q \sum_{j\in J_1'} m_j^{3/2}\,.
\]

Now, suppose that we've made $k$ steps of this construction. If
$J_k = \{1, \, ...\,, N\}$, then we are done with $J'=J_k'$. If
$\{1, \, ...\,, N \} \setminus J_k \ne \varnothing$, we choose
$j_{k+1}\in \{1, \, ...\,, N\}\setminus J_k $ such that
\[
m_{j_{k+1}} = \max \big\{ m_j\colon j\in \{1, \, ...\,,
N\}\setminus J_k \big\}\,,
\]
and define the sets $J_{k+1}' = J_k' \cup \{j_{k+1}\}$,
\[
J_{k+1}'' = J_k'' \cup \big\{j\in \{1, \, ...\, N\}\setminus J_k
\colon 0<|j-j_{k+1}|_* < 2Q\sqrt{m_{j_{k+1}}} \big\},
\]
and $J_{k+1} = J_{k+1}' \cup J_{k+1}''$. Then, as above,
\[
\sum_{j\in J_{k+1} \setminus J_k} m_j \le 5Q m_{j_{k+1}}^{3/2}\,,
\]
whence
\[
\sum_{j\in J_{k+1}} m_j \le 5Q \sum_{j\in J_{k+1}'} m_j^{3/2}\,.
\]
We are done. \hfill $\Box$

\subsection{Probabilistic preliminaries}\label{sect2.1}

\begin{lemma}{\rm \cite[Lemma~2.1]{NSV}}\label{lemma2.1}
Let $\eta_k$ be standard complex Gaussian random variables (not
necessarily independent). Let $a_k>0$, $ S=\sum_k{a_k}$. Then, for
every $t>0$,
\[
\Pr  { \sum_{k}a_k|\eta_k|>t }  \le 2e^{-\frac 12 (t/S)^2}\,.
\]
\end{lemma}

We also need the following  classical Bernstein's estimate:
\begin{lemma}\label{lemma2.3}
Let $\psi_k$, $k=1, 2, \, ... \,, n$,  be independent random
variables with zero mean such that, for some $K>0$ and every $t>0$,
\[
\Pr {|\psi_k|>t} \le K e^{-t}\,.
\]
Then, for $0 < t \le 5Kn$,
\[
\Pr { \big| \sum_k \psi_k \big| > t } \le 2 \exp\left(-
\frac{t^2}{16K n}\right).
\]
\end{lemma}

\par\noindent{\em Proof:} Set $\displaystyle S_n = \sum_{k=1}^n
\psi_k$.
Then $\displaystyle \Ex e^{\lambda S_n} = \prod_{k=1}^n \Ex
e^{\lambda \psi_k}$. Note that
\begin{multline*}
\Ex e^{\lambda \psi_k} = \Ex\left\{1 + \lambda \psi_k + \left(
e^{\lambda\psi_k}-1-\lambda\psi_k \right) \right\} \\
= 1 + \lambda \left\{ \int_{0}^{\infty} \Pr { \psi_k>t } \left(
e^{\lambda t}-1 \right)\,dt + \int_{0}^{\infty} \Pr {
\psi_k<-t} \left( 1-e^{-\lambda t} \right)\, dt \right\} \\
\le 1 + K\lambda\left\{ \int_0^\infty e^{-t} \left( e^{\lambda
t}-1 \right)\, dt + \int_0^\infty e^{-t} \left( 1-e^{-\lambda t}
\right)\, dt  \right\} \\
=  1 + \frac{2K\lambda^2}{1-\lambda^2} \le 1 + 4K \lambda^2 \le
e^{4K \lambda^2}\,,
\end{multline*}
provided that $\lambda \le \tfrac23$. Hence, we get
\[
\Pr { S_n > t } \le e^{-\lambda t} \Ex e^{\lambda S_n } \le e^{4K
n \lambda^2 - \lambda t}\,.
\]
Similarly, $ \Pr { S_n<-t } \le e^{4K n \lambda^2 - \lambda t}$,
and therefore $ \Pr { |S_n|>t } \le 2e^{4K n \lambda^2 - \lambda
t}$. Taking $\displaystyle \lambda = \frac{t}{8Kn}$, we get the
lemma. \hfill $\Box$

\subsection{Mean number of zeroes of a Gaussian Taylor series}
\label{sect2.1.1}

Consider a Gaussian Taylor series
\[
g(z) = \sum_{k=0}^\infty \zeta_k a_k z^k
\]
with non-negative $a_k$ such that $\displaystyle
\lim_{k\to\infty} \sqrt[k]{a_k} = 0$. Then almost surely, the
series on the right-hand side has infinite radius of convergence,
and hence $g$ is an entire function. By $n_g(r)$ we denote the
number of zeroes of the function $g$ in the disk of radius $r$.

\begin{lemma}\label{lemma2.EK}
\[
\Ex n_g(r) = \frac12\, \frac{r\mathfrak C_g'(r)}{\mathfrak
C_g(r)}\,,
\]
where
\[
\mathfrak C_g (r) = \sum_{k=0}^\infty a_k^2 r^{2k}\,.
\]
\end{lemma}

This readily follows from the Edelman-Kostlan formula for the
density of mean counting measure of zeroes of an arbitrary
Gaussian analytic function, see \cite[Section~2]{S}.
Alternatively, one can obtain this formula using the argument
principle, see \cite[page 195, Exercise~5]{Kahane}.

\subsection{Operators $T_w$ and shift invariance}\label{sect2.2}

For a function $g\colon \C\to\C$ and a complex number $w\in \C$,
we define
$$
T_w g(z)=g(w+z)e^{-z\overline w}e^{-\frac 12|w|^2}\,.
$$
In what follows, we use some simple properties of these operators.
\begin{itemize}
\item[(a)] $T_w$ are unitary operators in the Fock-Bargmann
space of entire functions:
\begin{multline*}
\| T_w f \|^2 = \frac1{\pi} \iint_{\C} |f(w+z)|^2 e^{-2{\rm Re}\,
(z\overline{w}) -|w|^2 - |z|^2}\, dm(z) \\
= \frac1{\pi} \iint_{\C} |f(w+z)|^2 e^{-|w+z|^2}\, dm(z) =
\|f\|^2\,.
\end{multline*}
\item[(b)] If $f$ is a G.E.F., then $T_w f$ is a G.E.F. as
well.
\end{itemize}
In particular, the distribution of the random zero set
$\Z_f=f^{-1}\{0\}$ is translation invariant. The property  (b)
also yields the distribution invariance of the function $f^*(z) =
|f(z)|e^{-|z|^2/2}$ with respect to the isometries of $\C$.
Indeed, a straightforward inspection shows that $(T_w f)^*(z) =
f^*(w+z)$.
\begin{itemize}
\item[(c)] By (b), if $f$ is a G.E.F., then
\[
T_w f = \sum_{k\ge 0} \zeta_k(w) \frac{z^k}{\sqrt{k!}} \,,
\]
where $ \zeta_k(w) $ are independent standard complex Gaussian
random variables. Recalling that $\displaystyle \Bigl\{
\frac{z^k}{\sqrt{k!}} \Bigr\}_{k\ge 0}$ is an orthonormal basis in
the Fock-Bargmann space, and using that $T_w$ is a unitary
operator and $T_w T_{-w}$ is the identity operator in that space,
we get
\[
\zeta_k (w) = \Bigl\langle T_w f, \frac{z^k}{\sqrt{k!}}
\Bigr\rangle = \Bigl\langle f, T_{-w}\Bigl( \frac{z^k}{\sqrt{k!}}
\Bigr) \Bigr\rangle.
\]
\end{itemize}
Note that for $w\ne w'$, the Gaussian variables $\zeta_k(w)$ and
$\zeta_{k'}(w')$ are correlated and
\begin{itemize}
\item[(d)]
\[
\Bigl| \Ex \big\{ \zeta_k(w) \overline{\zeta_{k'}(w')} \big\}
\Bigr| = \left| \Bigl\langle T_{-w}\Bigl( \frac{z^k}{\sqrt{k!}}
\Bigr), T_{-w'}\Bigl( \frac{z^{k'}}{\sqrt{k'!}} \Bigr)
\Bigr\rangle \right|\,.
\]
\end{itemize}

\bigskip Let $\gamma \subset \C$ be an oriented  curve. For any
continuous
function $h\colon \gamma \to \C\setminus \{0\}$, we denote by
$\Delta_{\gamma} \arg h$ the {\em increment} of the argument of
$h$ over $\gamma$. Note that if $f$ does not vanish on the curve
$\gamma$, then
\[
\Delta_{\gamma-w} \arg T_w f = \Delta_\gamma \arg f -
\Delta_\gamma {\rm Im} (z-w)\overline{w} = \Delta_\gamma \arg f -
\Delta_\gamma {\rm Im} (z\overline{w})\,.
\]
Here and in what follows, $\gamma-w$ denotes the translation of the
curve $\gamma$ by $-w$.
\begin{itemize}
\item[(e)] Set $\delta(f, \gamma) = \Delta_\gamma
\arg f - \Ex \Delta_\gamma \arg f $. Then $ \delta(T_w f,
\gamma-w) = \delta (f, \gamma)$.
\end{itemize}

\medskip If $I\subset R\T$ is an arc centered at $w$ and oriented
counterclockwise, then using rotation invariance and the argument
principle, we get \[ \Ex \Delta_I \arg f = \frac{|I|}{2\pi R}\, \Ex
\Delta_{R\T} \arg f = \frac{|I|}{2\pi R}\, \Ex 2\pi n(R) =
\frac{|I|}{R}\, R^2 = |I| R\] and \[ \delta (T_w f, I-w) = \Delta_I
\arg f - |I| R\,.\]

\subsection{Almost independence}\label{sect2.3}

Our approach is based on the almost independence property
introduced in \cite{NSV}. It says that if $\{w_j\}\subset \C $ is
a ``well-separated'' set, then the G.E.F. $T_{w_j}f$ can be
simultaneously approximated by independent G.E.F.. The following
lemma somewhat extends Theorem~3.2 from \cite{NSV}.

\begin{lemma}\label{thm3.1}
There exists a numerical constant $A>1$ such that for every family
of pairwise
disjoint disks $D(w_j, r_j+A\rho_j)$ with
\[
w_j\in\C, \quad r_j\ge 1, \quad \rho_j\ge \max\bigl( 1, \sqrt{\log
r_j}\, \bigr)\,,
\]
one can represent the family of G.E.F. $T_{w_j}f$ as
\[
T_{w_j} f = f_j + h_j
\]
where $f_j$ are independent G.E.F. and
\[
\Pr { \max_{z\in r_j\D} |h_j(z)| e^{-|z|^2/2} \ge e^{-\rho_j^2} }
\le 2\exp\bigl( -\tfrac12 e^{\rho_j^2} \bigr) \,.
\]
\end{lemma}
Theorem~3.2 in \cite{NSV} corresponds to the case when $r_j= r\ge 1$
and $\rho_j = Nr$ with $N\ge 1$. We prove Lemma~\ref{thm3.1} in the
Appendix.

\subsection{Bounds for G.E.F.}\label{sect2.4}
Our first lemma estimates the probability that the function $f$ is
very large:
\begin{lemma}{\rm (cf. \cite[Lemma~4.1]{NSV})}\label{cor2.2}
Let $f$ be a G.E.F.. Then, for each $r \ge 1$ and $M\ge 1$,
\[
\Pr { \max_{z\in r\D} |f(z)|e^{-|z|^2/2}  \ge M } \le 18r^2
e^{-\tfrac1{32} M^2}\,.
\]
\end{lemma}

\par\noindent{\em Proof:} We cover the disk $r\D$ by at most
$(2r+1)^2 \le 9r^2$ disks $\mathcal D_j$ of radius $1$ and show
that for each $j$,
\[
\Pr { \max_{z\in \mathcal D_j} |f(z)|e^{-|z|^2/2}  \ge M } \le
2e^{-\tfrac1{32} M^2}\,.
\]
By the translation invariance of the distribution of the random
function $|f(z)|e^{-|z|^2/2}$ it suffices to prove this estimate
in the unit disk $\D$. Clearly,
\begin{multline*}
\Pr { \max_{z\in \D} |f(z)|e^{-|z|^2/2}  \ge M } \le \Pr {
\max_{z\in \D} |f(z)| \ge M }
\\
\le \Pro { \sum_{k\ge 0} \frac{|\zeta_k|}{\sqrt{k!}} \ge M }
\stackrel{\rm Lemma~\ref{lemma2.1}}\le 2e^{-\tfrac12 (M/S)^2}
\end{multline*}
with $
S = \sum_{k\ge 0} \frac1{\sqrt{k!}} < 4$.
Hence, the lemma. \hfill $\Box$

\bigskip The following lemma estimates the probability that the
function
$f$ is very small:
\begin{lemma}{\rm (cf. Lemma~8 in \cite{Krishnapur} and
Lemma~4.2 in \cite{NSV})}\label{lemma4.1} Let $f$ be a G.E.F.. Let
$r\ge 1$ and $m\ge 3$. Then
\[
\Pr { \max_{r\D} |f| \le e^{-m r^2} } \le \exp \Bigl( -
\frac{m^2}{\log m}\, r^4 \Bigr)\,.
\]
\end{lemma}

\par\noindent{\em Proof:} Suppose that $|f|\le e^{-m r^2}$
everywhere in $r\D$. Then by Cauchy's inequalities,
\[
|\zeta_n| \le \frac{\sqrt{n!}}{r^n} \max_{r\D} |f| \le
\frac{n^{n/2}}{r^n} e^{-m r^2}, \qquad n = 0, 1, 2,\, ... \,.
\]
For $\displaystyle 0 \le n \le \frac{m}{\log m}\, r^2 $, the
probabilities of these events do not exceed
\[
\left( n r^{-2} \right)^n e^{-2m r^2} \le \left( \frac{m}{\log m}
\right)^{\frac{m}{\log m}\, r^2} e^{-2m r^2} < e^{-mr^2}.
\]
Since these events are independent, the probability we are
estimating is bounded by
\[
\exp\left( -m r^2 \Bigl( \frac{m}{\log m}\, r^2\Bigr) \right) =
\exp \Bigl( - \frac{m^2}{\log m} r^4 \Bigr)\,.
\]
We are done. \hfill $\Box$

\bigskip The next lemma  bounds the probability that
a G.E.F. is small on a given curve of a given length.
\begin{lemma}\label{lemma4.2}
Let $f$ be a G.E.F., and let $\gamma$ be a curve of length at most
$r\ge 1$. Then, for any positive $\epsilon\le \tfrac14$,
\[
\Pr { \min_{z\in \gamma} |f(z)|e^{-|z|^2/2} < \epsilon } < 100 r
\epsilon \sqrt{\log\frac1{\epsilon}}\,.
\]
\end{lemma}

\par\noindent{\em Proof:} We split the curve $\gamma$ into $\lceil
r \rceil$ arcs $\gamma_j$ of length at most $1$, and fix the
collection of disks $\mathcal D_j$ of radius $1$ such that
$\gamma_j \subset \mathcal D_j$. We'll show that for each $j$,
\[
\Pr { \min_{z\in\gamma_j} |f(z)|e^{-|z|^2/2} < \epsilon  } < 50
\epsilon \sqrt{\log\frac1{\epsilon}}\,.
\]
Clearly, this will yield the lemma.

By the shift invariance of the distribution of the random function
$|f(z)|e^{-|z|^2/2}$, we assume without loss of generality that
$\mathcal D_j $ is the unit disk $\D$. Taking into account that $
e^{-|z|^2/2} > \tfrac12 $ everywhere in the unit disk, we have
\[
\Pr { \min_{z\in\gamma_j} |f(z)|e^{-|z|^2/2} < \epsilon  } \le \Pr
{ \min_{z\in\gamma_j} |f(z)| < 2\epsilon }\,.
\]

We choose points $\{z_m\}\subset\gamma$ and disks $D_m =
\{|z-z_m|\le \kappa\epsilon\}$ such that
\[
\gamma \subset \bigcup_m D_m\,, \quad {\rm and} \quad \#\{z_m\}
\le \Big\lceil \frac1{2\kappa\epsilon} \Big\rceil\,,
\]
with the parameter $\kappa$ to be specified later. Then, for $z\in
D_m$,
\[
|f(z)| \ge |f(z_m)| - |z-z_m| \, \max_{\D}|f'| \ge |f(z_m)| -
\kappa\epsilon \max_\D |f'|\,.
\]
Hence, we need to estimate the probability of the events
\[
\Omega_1 = \Big\{  \min_m |f(z_m)|\le 3\epsilon \Big\} \quad {\rm
and } \quad
\Omega_2 = \Big\{ \max_\D |f'| \ge \frac1{\kappa}\Big\}\,.
\]
If neither of these events holds, then $|f(z)| > 3\epsilon -
\epsilon = \epsilon$ everywhere on $\gamma$.

Recall that for any standard complex Gaussian random variable
$\zeta$ and for any $t>0$, we have $\Pr {|\zeta|\le t} < t^2
$, also recall that $f(z_m)e^{-|z_m|^2/2}$ is a standard complex
Gaussian random variable. Hence, for any fixed $m$, we have $\Pr {
|f(z_m)| \le 3\epsilon} \le \Pr { |f(z_m)|e^{-|z_m|^2/2} \le
3\epsilon }  < 9\epsilon^2$. Therefore,
\[
\Pr { \Omega_1 } < \Big\lceil
\frac1{2\kappa \epsilon} \Big\rceil \cdot 9 \epsilon^2  \le
\frac92 \epsilon\kappa^{-1} + 9\epsilon^2\,.
\]

Next,
\[
\Pr { \Omega_2 } \le \Pr { \sum_{k\ge 1} \frac{k}{\sqrt{k!}}
|\zeta_k| \ge \frac1{\kappa} } \stackrel{\rm
Lemma~\ref{lemma2.1}}\le
2e^{-\tfrac12 (\kappa S)^{-2}}
\]
with  $S = \sum_{k\ge 1} \frac{k}{\sqrt k!} < 6$.
Therefore, $
\Pr { \Omega_2 } \le 2e^{-\tfrac1{72} \kappa^{-2}}$,
and
\[
\Pr {\Omega_1} + \Pr {\Omega_2} < \frac92 \epsilon \kappa^{-1} +
2e^{-\tfrac1{72} \kappa^{-2}} + 9\epsilon^2\,.
\]
Choosing here $ \kappa^{-1} = \sqrt{72\log\frac1{\epsilon}}$,
we get
\begin{multline*}
\Pr { \min_{z\in\gamma_j} |f(z)| < 2\epsilon } \le \Pr {\Omega_1}
+ \Pr {\Omega_2} \\
< 27\sqrt{2}\, \epsilon\sqrt{\log\frac1{\epsilon}} + 2\epsilon +
9\epsilon^2
< 50 \epsilon\sqrt{\log\frac1{\epsilon}}\,,
\end{multline*}
proving the lemma. \hfill $\Box$

\subsection{Upper bounds for the increment of the
argument}\label{sect2.6}

We say that a piecewise $C^1$-curve $\gamma\subset r\D$ is {\em
good} if its length does not exceed $r$ and, for any $\zeta\in\C
\setminus \{\gamma\}$, we have $\big| \Delta_\gamma \arg (z-\zeta)
\big| \le 2\pi$. The following lemma is classical (cf.
\cite[Lemma~6, Chapter~VI]{BYa}):

\begin{lemma}\label{lemma6.2} There exists a numerical
constant $B>1$ with the following property. Let $g$ be an analytic
function in the disk $2r\D$ such that $\displaystyle \sup_{2r\D}
|g|
\le 1$.
If $\displaystyle \max_{r\D} |g| \ge e^{-\beta}$, then for any good
curve $\gamma\subset r\D$, we have
\[
\big| \Delta_\gamma \arg g \big| \le B \beta\,.
\]
\end{lemma}

\par\noindent{\em Proof:} By scale invariance, it suffices to prove
the lemma for $r=1$. Choose $z_0\in r\T$ such that $\displaystyle |g(z_0)| = \max_{r\D} |g| \ge e^{-\beta}$, and denote by $\phi$ a M\"{o}bius transformation $\phi\colon 2\D\to 2\D$ with $\phi (0) = z_0$. By Jensen's formula applied to the function $g \circ \phi$, the number of zeroes of $g$ in the disk $\tfrac32 \D$ does not exceed $C_1 \beta$. Hence,
$g=pg_1$ where $p$ is a polynomial of degree $N\le C_1\beta$ with
zeroes in $\tfrac32 \D$ and a unimodular leading coefficient, and
$g_1$ does not vanish in $\tfrac32 \D$, $g_1(0)>0$.

\medskip\par\noindent{\bf Claim~9\,-1.} $\displaystyle
\int_0^{2\pi}
\big| \log|g_1(\tfrac32 e^{i\theta})| \big|\, \frac{d\theta}{2\pi}
\le C_2 \beta $.

\medskip\par\noindent{\em Proof of Claim~9\,-1:} Indeed,
\begin{multline*}
\int_0^{2\pi} \big| \log|g_1(\tfrac32 e^{i\theta})| \big| \,
\frac{d\theta}{2\pi} \le \int_0^{2\pi} \big| \log|g(\tfrac32
e^{i\theta})| \big| \, \frac{d\theta}{2\pi} + \int_0^{2\pi} \big|
\log|p(\tfrac32 e^{i\theta})| \big| \, \frac{d\theta}{2\pi}
\\
\le \int_0^{2\pi} \log^-|g(\tfrac32 e^{i\theta})|  \,
\frac{d\theta}{2\pi} + C_3 \beta\,.
\end{multline*}
To estimate the integral on the right-hand side, we note that
\[
\int_0^{2\pi} \log |g(\tfrac32 e^{i\theta})|\, \frac{\bigl(
\tfrac32 \bigr)^2- |z_0|^2}{|\frac32 e^{i\theta} - z_0|^2} \,
\frac{d\theta}{2\pi} \ge \log |g(z_0)| \ge -\beta\,,
\]
whence
\[
\int_0^{2\pi} \log^- |g(\tfrac32 e^{i\theta})|\,
\frac{d\theta}{2\pi} \le C_4\beta\,.
\]
Hence, the claim. \hfill $\Box$

\bigskip Now, $\big| \Delta_\gamma \arg g \big| \le \big|
\Delta_\gamma \arg p \big| + \big| \Delta_\gamma \arg g_1 \big| $,
and since the curve $\gamma$ is good, we have $ \big|
\Delta_\gamma \arg p \big| \le 2\pi N \le C_5\beta$. Consider the
function $h=\arg g_1$ normalized by the condition $h(0)=0$. The
function $h$ is harmonic in $\tfrac 32 \D$. We have
\[
\big | \Delta_\gamma h \big| \le \max_{\overline{\D}} h -
\min_{\overline{\D}} h \le 2 \max_{\overline{\D}} |h|\,.
\]
Since
\[
h(z) = \int_0^{2\pi} \log |g_1(\tfrac32 e^{i\theta})| \, {\rm Im}
\, \frac{\tfrac32 e^{i\theta} + z}{\tfrac32 e^{i\theta} - z}\,
\frac{d\theta}{2\pi}\,, \qquad |z|\le 1\,,
\]
we have
\[
|h(z)| \le C_6 \int_0^{2\pi} \big| \log |g_1(\tfrac32 e^{i\theta}
)|\, \big|\, d\theta \stackrel{\rm Claim~9-1}\le C_7\beta\,,
\qquad
|z|\le 1\,.
\]
This proves the lemma. \hfill $\Box$

\begin{lemma}\label{lemma.new-old}
Let $r\ge 1$, let $\gamma\subset r\D$ be a good curve, let
$m\ge 25 B$, and let $f$ be a G.E.F.. Consider the event $\Omega =
\bigl\{ | \delta (f, \gamma)| \ge mr^2 \bigr\} $. Then
\[
\Omega \subset \Omega' \cup \Bigl\{ \max_{r\D} |f| < e^{-\frac1{4B}
mr^2 } \Bigr\} \qquad {\rm with} \quad
\Pr {\Omega'} \le \exp\bigl(-e^{\frac1{6B} mr^2} \bigr)\,.
\]
In particular,
\[
\Pr {\Omega } \le 2\exp\bigl( -\tfrac1{16B^2} \,
\tfrac{m^2r^4}{\log m} \bigr)\,.
\]
\end{lemma}

\par\noindent{\em Proof:} Introduce the events
\[
\Omega_1(m) = \bigl\{ | \Delta_\gamma \arg f | \ge mr^2 \bigr\}\,,
\]
and
\[
\Omega' (m) = \bigl\{ \max_{z\in 2r\D} |f(z)| e^{-|z|^2/2}>e^{
\frac1{3B} mr^2} \bigr\}\,.
\]

\medskip\par\noindent{\bf Claim 10\,-1.} {\em For $m\ge 12B$,}
$\Omega_1 (m) \subset \Omega'(m) \cup \bigl\{ \max_{r\D} |f| <
e^{-\frac1{2B} mr^2}\bigr\} $.

\medskip\par\noindent{\em Proof of Claim~10\,-1:} Suppose that the
event $\Omega'(m)$ does not occur. Then
\[
\max_{2r\D}|f| \le e^{\frac1{3B} mr^2 + 2r^2 } =
e^{( \frac1{3B} + \frac2{m}) mr^2}
\stackrel{m\ge 12B} \le e^{\frac1{2B} mr^2}\,.
\]
If the event $\Omega_1 (m)$ occurs, then by Lemma~\ref{lemma6.2}
\[
mr^2 \le \bigl| \Delta_\gamma \arg f \bigr| \le B\,
\log\frac{\max_{2r\D}|f|}{\max_{r\D}|f|}\,,
\]
whence,
\[
\max_{r\D} |f| \le e^{-\frac1{B} mr^2} \max_{2r\D} |f| \le
e^{-\frac1{2B} mr^2}\,,
\]
proving the claim. \hfill $\Box$

\medskip\par\noindent{\bf Claim 10\,-2. }{\em For $m\ge 12B$,}
$ \Pr { \Omega'(m) } \le \exp\bigl( -e^{\frac1{3B} mr^2} \bigr)$.

\medskip\par\noindent{\em Proof of Claim~10\,-2:} We have
\[
\Pr {\Omega'(m)} \stackrel{\rm Lemma~\ref{cor2.2}}\le 72r^2
\exp\Bigl( -\tfrac1{32}  e^{\frac2{3B} mr^2 } \Bigr)
\le \tfrac6{B} mr^2
\exp\Bigl( -\tfrac1{32} e^{\frac2{3B} mr^2} \Bigr)\,.
\]
It's easy to see that for $t=\tfrac1{3B} mr^2 \ge 4$, one has
\begin{multline*}
18t \exp\bigl( -\tfrac1{32} e^{2t} \bigr) \le 18t \exp\bigl(
-\tfrac{e^4}{32} e^t \bigr) < 18t \exp\bigl( -\tfrac32 e^t \bigr) \\
= \underbrace{18t \exp\bigl( -\tfrac12 e^t \bigr)}_{<1}\,  \exp\bigl(
- e^t \bigr) < \exp\bigl( -e^t \bigr)\,.
\end{multline*}
Hence, the claim. \hfill $\Box$

\medskip\par\noindent{\bf Claim 10\,-3. } {\em For $m\ge 12B$, }
$\Pr {\Omega_1 (m)} \le 2 \exp\bigl( -\tfrac1{4B^2}\,
\frac{(mr^2)^2}{\log (m/2B)}\bigr) $.

\medskip\par\noindent{\em Proof of Claim~10\,-3:} By
Lemma~\ref{lemma4.1},
\[
\Pr { \max_{r\D} |f| < e^{-\frac1{2B} mr^2} }
\le \exp\bigl( -\tfrac1{4B^2}
\tfrac{(mr^2)^2}{\log(m/2B)} \bigr)\,.
\]
It's easy to check that for $t=\tfrac1{B} mr^2\ge 12$, one has
$e^{t/3} > \tfrac{t^2}4$. Therefore,
\[
e^{\frac1{3B} mr^2} > \tfrac1{4B^2} (mr^2)^2
> \tfrac1{4B^2}\, \tfrac{(mr^2)^2}{\log(m/2B)}\,.
\]
Thus $ \Pr { \Omega'(m) } $ also does not exceed $\exp\bigl(
-\tfrac1{4B^2} \tfrac{(mr^2)^2}{\log(m/2B)} \bigr) $. We are done.
\hfill $\Box$

\medskip\par\noindent{\bf Claim 10\,-4.} $\bigl| \Ex \Delta_\gamma
\arg f \bigr| \le 12.5 Br^2$.

\medskip\par\noindent{\em Proof of Claim~10\,-4:} By Claim~10\,-3,
for $s\ge 12Br^2$, we have
\[
\Pr { |\Delta_\gamma \arg f| \ge s } \le  2 \exp\Bigl(
\frac{(s/2B)^2}{\log s/(2Br^2) }\Bigr) \le 2 \exp\Bigl(
\frac{(s/2B)^2}{\log s/2B }\Bigr).
\]
Therefore,
\begin{multline*}
\bigl| \Ex \Delta_\gamma \arg f \bigr| \le 12 Br^2 + 2
\int_{12Br^2}^\infty \exp\Bigl( - \frac{(s/2B)^2}{\log
(s/2B)}\Bigr)\, ds \\
= 12Br^2 + 4B \underbrace{\int_{6r^2}^\infty e^{-s^2/\log s}\,
ds}_{<1/8} < 12.5Br^2\,,
\end{multline*}
proving the claim. \hfill $\Box$

\medskip Now, we readily finish the proof of
Lemma~\ref{lemma.new-old}. Suppose that the event $\Omega$ occurs;
i.e., $\bigl| \delta(f, \gamma) \bigr| \ge mr^2$ with $m\ge 25B$.
Then
\[
\bigl| \Delta_\gamma \arg f \bigr| \ge  \bigl| \delta(f, \gamma)
\bigr| -  \bigl| \Ex \Delta_\gamma \arg f \bigr| \ge (m-12.5B)r^2
\ge
\tfrac12 mr^2\,.
\]
That is, $\Omega \subset \Omega_1 (\tfrac12 m)$ and
the lemma follows from Claims~10\,-1, 10\,-2 and 10\,-3 applied
with
$\tfrac12 m$ instead of $m$. \hfill $\Box$

\medskip\par\noindent{\em Remark.} If $\gamma\colon [0, 1]\to \C$ is a good curve, then
\[
\Ex \Delta_\gamma \arg f = \operatorname{Im}\, \int_0^1 \gamma'(t) \overline{\gamma (t)}\, dt\,.
\]
Then taking into account that $\Ex \frac{f'}{f}(z) = \overline{z}$, one can get a better estimate $\bigl| \Ex \Delta_\gamma \arg f \bigr| \le r^2$ than the one given in Claim 10\,-4.

\section{The upper bound for $1<\alpha<2$}\label{sect3}

\subsection{Few arcs with large increments of the
argument}\label{subsect3.1}

Given $r\ge 1$, we fix a collection of $N\le 2\pi \frac{R}{r}$
disjoint arcs $\big\{I_j\big\}_{1\le j \le N}$ of length $r$ on
the circumference $R\T$. Then, given $\Lambda \ge 1$ and a
positive integer $L$, we introduce two events. The first event
$\Omega_1 (r, R, \Lambda, L)$ is that the collection
$\big\{I_j\big\}_{1\le j \le N}$ contains a sub-collection of $L$
disjoint arcs $\{I_j\}_{j\in J}$ such that
\[
\sum_{j\in J} \big| \delta(f, I_j) \big| \ge \Lambda\,.
\]
To define the second event, we fix $N$ independent G.E.F. $f_j$.
Then the event $\Omega_2 (r, R, \Lambda, L) $ is that the
collection $\big\{I_j\big\}_{1\le j \le N}$ contains a
sub-collection of $L$ disjoint arcs $\{I_j\}_{j\in J}$ such that,
\[
\sum_{j\in J} \big| \delta (f_j, \widetilde{I_j}) \big| \ge
\Lambda\,.
\]
Here, $\widetilde{I_j} = I_j - w_j$, where $w_j$ are the centers of
the arcs $I_j$.

\begin{lemma}\label{prop7.1} Suppose that $R$ is sufficiently big.
Suppose also that
\[
R^{1/2} \le \Lambda \le R^2 \qquad {\rm  and} \qquad 1 < L <
\frac{b \Lambda}{r^2 + \log R}
\]
with a sufficiently small positive numerical constant $b$. Then
the probabilities of the events $\Omega_i$, $i=1, 2$, do not
exceed $e^{-b_1r^2 \Lambda}$ with a positive numerical constant
$b_1$.
\end{lemma}

\par\noindent{\em Proof of Lemma~\ref{prop7.1}}: First, we estimate
the
probability of the event $\Omega_2$; this is a simpler part of
the job. Suppose that the event $\Omega_2 (r, R, \Lambda, L)$
occurs.
We choose $ M_j \le r^{-2}\, \big| \delta (f_j, \widetilde{I}_j)
\big|$
such that $ \sum_{j\in J} M_j r^2 = \Lambda$. Let $B$ be the
constant from Lemma~\ref{lemma6.2}.
Note that the arcs $\bigl\{ \widetilde{I}_j \bigr\}$ with $M_j<50B$
can contribute at most
$ 50BL<50Bb\Lambda<\tfrac12 \Lambda $ to the total sum, provided
that $b<\tfrac1{100B}$.
We discard the arcs $\bigl\{ \widetilde{I}_j \bigr\}$ with
$M_j<50B$ and denote by $J$ the collection of remaining arcs.

Now, let $m_j$ be the largest positive integer power of $2$ such
that $m_j\le M_j$, $j\in J$. Then
\begin{equation}\label{eq.verynew}
\frac14 \Lambda \le \sum_{j\in J} m_j r^2 \le \Lambda \qquad {\rm
and} \qquad m_j\ge 25B\,,
\end{equation}
and
\[
\Pr { \Omega_2 } \le \sum_{J} \sum_{\{m_j\}}
\Pro { \bigcap_{j\in J} \bigl\{ |\delta(f_j, \widetilde{I}_j)| \ge
m_j r^2 \bigr\}  }
\]
where the first sum is taken over all subsets $J\subset \{1, ... ,
N\}$ of cardinality
at most $L$, and the second sum is taken over all possible choices
of $m_j$, $j\in J$,
that are positive integer powers of $2$ satisfying restrictions
\eqref{eq.verynew}.
Since $f_j$ are independent, we have
\[
\Pro { \bigcap_{j\in J} \bigl\{ |\delta(f_j, \widetilde{I}_j)| \ge
m_j r^2 \bigr\}  }
= \prod_{j\in J} \Pr { |\delta(f_j, \widetilde{I}_j)| \ge m_j r^2
}\,.
\]
The probabilities of the events on the right-hand side were
estimated in Lemma~\ref{lemma.new-old}:
\[
\Pr {   |\delta(f_j, \widetilde{I}_j)| \ge  m_j r^2 } \le
2\exp \Bigl( -\frac1{16B^2}\, \frac{m_j^2 r^4}{\log m_j}  \Bigr)\,.
\]
Therefore,
\begin{multline*}
\Pro { \bigcap_{j\in J} \bigl\{ |\delta(f_j, \widetilde{I}_j)| \ge
m_j r^2 \bigr\}  }
\le 2^L \exp\Bigl( -\frac1{16B^2}\, r^2 \sum_{j\in J} \frac{m_j^2
r^2}{\log m_j}  \Bigr)
\\ < 2^L \exp\Bigl( -\frac1{16B^2}\, r^2\sum_{j\in J} m_j r^2
\Bigr)
\le
2^L \exp\Bigl( -\frac1{64B^2}\, r^2\Lambda\Bigr)\,,
\end{multline*}
and
\[
\Pr {\Omega_2 } < 2^L \exp\Bigl( -\frac1{64B^2}\, r^2\Lambda\Bigr)
\sum_J \sum_{\{m_j\}} 1\,.
\]

To get rid of the sums on the right-hand side, we need
to estimate the number of different ways to choose the ``data'' $J,
\{m_j\}_{j\in J}$\,.
Since $m_j$ is an integer power of $2$ and $m_j\le \Lambda$, for
each $j\in J$,
there are at most $2\log\Lambda$ ways to choose the integer $m_j$.
Hence, given set $J$ of cardinality at most $L$, we have at most
$(2\log\Lambda)^L$ ways to choose the collection $\{m_j\}_{j\in
J}$.
Also there are at most
\[
\sum_{0\le \ell \le L}\binom{N}{\ell} < (N+1)^{L} \stackrel{N\le
2\pi R}< e^{C L\log R}
\]
ways to choose the subset
$J\subset \{1, 2, \, ... \,, N\}$ of cardinality at most $L$.
Therefore,
\begin{eqnarray*}
2^L \sum_J \sum_{\{m_j\}} 1 &\le& (4\log\Lambda)^L e^{CL\log R} \\
&<& e^{C'L(\log R + \log\log\Lambda)} < e^{C''L\log R} <
e^{C''b\Lambda}\,,
\end{eqnarray*}
which is a negligible factor with respect to
$\exp\bigl( -\tfrac1{64B^2}\, r^2\Lambda\bigr)$, provided that
$b\ll B^{-2}$. This completes the estimate of $\Pr {\Omega_2}$.
\hfill $\Box$

\bigskip
The estimate of the probability of the event $\Omega_1$ follows a
similar pattern.
Now, the events  $\bigl\{ |\delta(f, I_j)| \ge m_j r^2 \bigr\}$,
$j\in J$,
are not independent. To get around this obstacle, we'll use the
almost independence lemma, which brings in some awkward
technicalities. We split the proof into several steps.

\medskip\par\noindent{\bf (i)} Suppose that the event
$\Omega_1 (r, R, \Lambda, L)$ occurs.
As above, we choose $ M_j \le r^{-2}\, \big| \delta (f, I_j) \big|$
such that $\sum_{j\in J} M_j r^2 = \Lambda$.
Then we fix a sufficiently large positive numerical constant $C_1
\ge 25B $ and note that the arcs $I_j$ with $M_j < 2C_1 (1 +
r^{-2}\log\Lambda)$ can contribute to the total deviation
$\Lambda$ at most
\[
2C_1 L ( r^2 + \log \Lambda ) < 2C_1\, \frac{b\Lambda (r^2 + 2\log
R)}{r^2 + \log R} \le 4b C_1 \Lambda
\]
which is much smaller than $\Lambda$ provided that the constant
$b$ is sufficiently small. We choose $b<\tfrac1{8C_1}$ and
conclude that {\em at least half of the deviation $\Lambda$ must
come from the arcs $I_j$ with sufficiently large $M_j$}. From now
on, we discard the arcs $I_j$ with $M_j < 2C_1 (1 +
r^{-2}\log\Lambda)$ and denote by $J$ the set of the remaining
arcs.

Now, let $m_j$ be the largest positive integer power of $2$ such
that $m_j\le M_j$, $j\in J$. Then
\begin{equation}\label{eq.cond}
\frac14 \Lambda \le \sum_{j\in J} m_j r^2 \le \Lambda\,, \qquad
{\rm
and} \qquad
m_j r^2 \ge C_1 \left( r^2 + \log\Lambda \right)
\end{equation}
and
\begin{equation}\label{eq.prob}
\Pr { \Omega_1 } \le \sum_{J} \sum_{\{m_j\}}
\Pro { \bigcap_{j\in J} \bigl\{ |\delta(f, I_j)| \ge m_j r^2
\bigr\}
}
\end{equation}
where the first sum is taken over all subsets $J\subset \{1, ... ,
N\}$ of cardinality
at most $L$, and the second sum is taken over all possible choices
of $m_j$, $j\in J$,
that are positive integer powers of $2$ satisfying restrictions
\eqref{eq.cond}.

As in the previous case, it suffices to show that, for a {\em
fixed}
subset $J \subset \{1, 2, \, ...\,, N\}$ with $\# J \le L$, and for
{\em fixed} $m_j$, $j\in J$, that are
integer powers of $2$ and satisfy conditions~\eqref{eq.cond}, one
has
\begin{equation}\label{eq3.3}
\Pro { \bigcap_{j\in J} \bigl\{ |\delta(f, I_j)| \ge m_j r^2
\bigr\}
} \le e^{-c
r^2 \Lambda}\,.
\end{equation}
Since we have at most $e^{C''L \log R} < e^{C''b\Lambda}$ possible
combinations of
the ``data'' $J$ and $\{m_j\}_{j\in J}$, the two sums
on the right-hand side of \eqref{eq.prob} contribute by
a negligible factor with respect to $e^{-c r^2 \Lambda}$,
provided that $b < \tfrac{c}{2C''}$.

\medskip\par\noindent{\bf (ii)} From now on, we fix  a
set $J$
of cardinality at most $L$, and $m_j$, $j\in J$, that are
integer powers of $2$ and satisfy conditions~\eqref{eq.cond}.
Let $w_j$ be the centers of
the arcs $I_j$, ${\widetilde I}_j = I_j-w_j$, and let
\[
\Omega_j \stackrel{\rm def} = \bigl\{ | \delta(f, I_j) | \ge m_jr^2
\bigr\} = \bigl\{ | \delta(T_{w_j}f, \widetilde{I}_j) | \ge m_jr^2
\bigr\}\,.
\]
By Lemma~\ref{lemma.new-old} applied to the G.E.F. $T_{w_j}f$ with
$\gamma=\widetilde{I}_j$ and $m=m_j$, we have
\[
\Omega_j \subset \Omega_j' \cup \Bigl\{ \max_{r\D} |T_{w_j}f| <
e^{-\frac1{4B} m_jr^2} \Bigr\}
\]
with $\Pr {\Omega_j'} < \exp \bigl( - e^{-\tfrac1{6B} m_j
r^2}\bigr)$, whence,
\[
\bigcap_{j\in J} \Bigl\{ \bigl| \delta(f, I_j) \bigr| \ge m_j r^2
\Bigr\} \subset \Bigl( \bigcup_{j\in J} \Omega_j' \Bigr) \cup
\bigcap_{j\in J} \Bigl\{ \max_{r\D} |T_{w_j}f| < e^{-\frac1{4B}
m_jr^2} \Bigr\}
\]
with
\[
\Pro { \bigcup_{j\in J} \Omega_j' } \stackrel{m_jr^2\ge
C_1\log\Lambda}\le  L \exp\Bigl( -e^{\frac1{6B} C_1\log\Lambda
}\Bigr) \stackrel{C_1\ge 25B}\le Le^{-\Lambda^4}
\stackrel{L<\Lambda}< \Lambda e^{-\Lambda^4}\,.
\]
Since $r^2 < \Lambda$, this is much less than $e^{-r^2\Lambda}$
when
$R\gg 1$.

Discarding the event $\displaystyle \bigcup_{j\in J} \Omega_j'$, we
need to estimate
the probability of the event
\[
\bigcap_{j\in J} \Bigl\{ \max_{r\D} |T_{w_j}f| \le \exp\bigl(
-\tfrac1{4B} m_jr^2\bigr) \Bigr\}
\]

\medskip\par\noindent{\bf (iii)} Combinatorial
lemma~\ref{claim7.2}
applied with $m_j=0$ for $j\notin J$ and with the constant
$Q=\frac{\pi}2 (A+1)$, gives us
a subset $J'\subset J$ such that
\[
|j-k|_* \ge Q(\sqrt{m_j} + \sqrt{m_k})\,, \qquad j,k\in J'\,, \
j\ne
k\,,
\]
and
\begin{equation}\label{eq.new}
\sum_{j\in J'} m_j^{3/2} \ge \frac1{5Q} \sum_{j\in J} m_j\,.
\end{equation}
Hence, the centers $w_j$ of the arcs from $J'$ are
well-separated:
\[
| w_j - w_k | = 2R \sin\frac{|j-k|_* r}{2R} \ge \frac2{\pi}\,
|j-k|_*\, r  \ge (A+1) \left( \sqrt{m_j}\, r + \sqrt{m_k}\, r
\right)
\]
for $j, k \in J'$, $j\ne k$. By the almost independence
lemma~\ref{thm3.1} applied with $r_j = \rho_j = \sqrt{m_j}r$, we
have $T_{w_j} f = f_j + h_j $, $j\in J'$, where $f_j$
are {\em independent} G.E.F., and
\begin{multline*}
\Pr { \max_{z\in r \D}|h_j(z)|e^{-|z|^2/2} \ge e^{-m_jr^2} } \le 2
\exp \big( -\tfrac12 e^{m_jr^2} \big) \\
\stackrel{m_jr^2\ge
C_1\log\Lambda}\le 2\exp\bigl( -\tfrac12 \Lambda^{C_1}\bigr)\,.
\end{multline*}
Introduce the event
\[
\mathcal F = \bigcup_{j\in J'} \Bigl\{  \max_{r\D} |h_j| >
\exp\bigl( -\tfrac12 m_j r^2\bigr) \Bigr\}\,.
\]

\medskip\par\noindent{\bf Claim~11\,-1. }
{\em For $R\gg 1$, $\Pr {\mathcal F} \le e^{-r^2\Lambda}$.}

\medskip\par\noindent{\em Proof of Claim~11\,-1:}
If for some $j\in J'$,
\[
\max_{r\D} |h_j| > \exp\bigl( -\tfrac12 m_j r^2\bigr)\,,
\]
then
\[
\max_{z\in r\D} |h_j(z)|e^{-|z|^2/2} > \exp\bigl( -\tfrac12
(m_j+1) r^2 \bigr) \stackrel{m_j \ge 25}> e^{-m_jr^2}.
\]
Therefore,
\begin{multline*}
\Pr {\mathcal F} \le \sum_{j\in J'} \Pr { \max_{z\in r\D}
|h_j(z)|e^{-|z|^2/2} > \exp\bigl( -\tfrac12 m_j r^2\bigr) } \\
\le L \cdot 2 \exp \big( -\tfrac12 \Lambda^{C_1} \big)
\stackrel{L\le\Lambda}\le 2\Lambda
\exp \big( -\tfrac12 \Lambda^{C_1} \big)\,.
\end{multline*}
Since $r^2<\Lambda$, this is much less than $e^{-r^2\Lambda}$,
provided that $R\gg 1$.
\hfill $\Box$

\medskip If the event $\mathcal F$ does not occur, then for each
$j\in J'$,
\begin{multline*}
\max_{r\D} |f_j| \le \max_{r\D} |T_{w_j} f| + \max_{r_j \D} |h_j|
\\
\le e^{ - \frac1{4B} m_j r^2} + e^{ - \frac12 m_jr^2}
\stackrel{B>1}< 2e^{-\frac1{4B} m_jr^2} \stackrel{m_jr^2\ge 25B}<
e^{ -\frac1{6B} m_j r^2}\,.
\end{multline*}
We conclude that {\em if $R$ is sufficiently big, then outside of
an
event of probability less than $\exp(-r^2\Lambda)$, we have
\[
\max_{r\D} |f_j| < \exp\bigl( -\tfrac1{6B}m_j r^2\bigr)
\]
for each $j\in J'$.}

\medskip\par\noindent{\bf (iv)} Our problem boils down to
the estimate of the probability that the {\em independent} events
\[
\Bigl\{ \max_{r\D} |f_j| < e^{- \frac1{6B} m_j r^2} \Bigr\}\,,
\qquad j\in J'\,,
\]
occur. By Lemma~\ref{lemma4.1} the logarithm of the probability of
each of these events doesn't exceed $ - \tfrac{c_2 m_j^2}{\log
m_j}\, r^4$ with $c_2 = \tfrac1{36B^2}$. Therefore, the logarithm
of the probability that all these events happen doesn't exceed
\begin{multline*}
-c_2 \sum_{j\in J'} \frac{m_j^2}{\log m_j}\, r^4 < -c_2 r^2
\sum_{j\in J'} m_j^{3/2} r^2 \\
\stackrel{\eqref{eq.new}}< - \frac{c_2}{5 \frac{\pi}2 (A+1)}\, r^2
\sum_{j\in J} m_j r^2 \stackrel{\eqref{eq.cond}}\le -c r^2
\Lambda
\end{multline*}
with
\[
c = \frac14\, \frac{2c_2}{5\pi (A+1)} = \frac1{360\pi B^2
(A+1)}\,.
\]
This completes the proof of \eqref{eq3.3} and, thereby, of the
lemma. \hfill $\Box$

\subsection{Proof of Theorem~\ref{thm_main}: the upper bound for
$1<\alpha<2$}
We need to estimate the probability of the event $ \Omega_\alpha =
\Bigl\{ |n(R)-R^2|>R^\alpha\Bigr\}$. Let $b$ be the constant from
the previous lemma. We fix a small positive $\delta\in (\tfrac14
b, \tfrac12 b)$ such that the number $N= \tfrac{2\pi}{\delta}
R^{2-\alpha}$ is an integer, take $r=\delta R^{\alpha-1}$, and
split the circumference $R\T$ into $N$ disjoint arcs $\{I_j\}$ of
length $r$. By the argument principle, $\Omega_\alpha = \Big\{
\big| \sum_j \delta (f, I_j) \big|
> 2\pi R^\alpha \Big\}$.
In the case $1<\alpha<2$, the cancelations between different random
variables $\delta
(f, I_j)$ are not important, so we are after the upper bound for
the probability of the bigger event $\Omega_\alpha' = \Bigl\{
\sum_j
\bigl| \delta (f, I_j) \bigr| > 2\pi R^\alpha \Bigr\}$.

We take $\Lambda = 2\pi R^\alpha$, and check that
Lemma~\ref{prop7.1} can be applied to the whole collection of arcs
$\{I_j\}$; i.e., with $L=N$. If $R$ is big enough then $\log R \ll
r^2$, and
\[
1 \le L = \frac{2\pi}{\delta} R^{2-\alpha} <  \frac12 \frac{b \cdot
2\pi
R^\alpha}{\delta^2 R^{2\alpha-2}} = \frac12 \frac{b\Lambda}{r^2} <
\frac{b\Lambda}{r^2 + \log R}\,.
\]
Therefore, the assumptions of Lemma~\ref{prop7.1} are fulfilled,
and we get
\[
\Pr { \Omega_\alpha' }
\le e^{- 2\pi b_1 r^2 R^\alpha} < e^{- c R^{3\alpha - 2 } }\,.
\]
Done! \hfill $\Box$

\section{The upper bound for $\tfrac12 < \alpha < 1$}\label{sect4}

\subsection{Approximating the total increment of $\arg f$
by the sum of increments of arguments of independent
G.E.F.}\label{subsect4.1}

\begin{lemma}\label{sub-thm2} Suppose that $R$ is sufficiently
big, that $1\le r\le 2$, and that $3R^{1/2}\le\Lambda \le
R$. Then, given a collection of disjoint arcs $\big\{ I_j \big\}$
of length $r$ of the circumference $R\T$ that are separated by
arcs of length at least $\log R$, there exists a collection of
independent G.E.F. $\{f_j\}$ such that
\[
\Pr { \big| \sum_{j} \delta(f, I_j) - \sum_j \delta(f_j,
\widetilde{I_j})  \big|  \ge \Lambda } \le e^{-b_2 \Lambda }\,,
\]
where $\widetilde{I_j} = I_j-w_j$ and $b_2$ is a positive numerical
constant.
\end{lemma}

\par\noindent{\em Proof of Lemma~\ref{sub-thm2}}: Set $\rho =
\sqrt{C_1\log R}$ with $C_1 \gg 1$. Let $A$ be the constant from
the
almost
independence lemma. If $R$ is big enough, then by
our assumptions, the disks $ D(w_j, r+A\rho)  $ are disjoint. So
the almost independence lemma~\ref{thm3.1} yields a decomposition
$T_{w_j}f = f_j + h_j $ with independent G.E.F. $\{f_j\}$ and
\[
\Pr { \max_{r\D} |h_j(z)|e^{-|z|^2/2}  \ge R^{-C_1} } \le 2 \exp
\big( -\tfrac12 R^{C_1} \big)\,.
\]
In what follows, we assume that
\[
\max_j \max_{r\D} |h_j(z)|e^{-|z|^2/2}  \le R^{- C_1}\,.
\]
For this, we throw away an event of probability at most
\[
2\pi R \cdot 2 e^{-\tfrac12 R^{C_1} } \ll e^{- \Lambda}\,.
\]

Since $\delta (f, I_j) = \delta(T_{w_j}f, \widetilde{I_j})$, we
need to estimate the probability of the event
\[
\Big\{ \Big| \sum_j \big[ \delta (T_{w_j}f, \widetilde{I_j}) -
\delta (f_j, \widetilde{I_j}) \big] \Big| \ge \Lambda \Big\}\,.
\]
Introduce the events
\[
\Omega_j = \big\{ \min_{z\in\widetilde{I_j}} |f_j(z)|e^{-|z|^2/2}
\le R^{-C_1/2} \big\}\,,
\]
and note that if $\Omega_j$ does not occur, then
\begin{multline*}
\big| \delta (T_{w_j}f, \widetilde{I_j}) - \delta (f_j,
\widetilde{I_j}) \big| = \big| \Delta_{\widetilde{I_j}} \arg
T_{w_j}
f - \Delta_{\widetilde{I_j}} \arg f_j \big| \\
= \big| \Delta_{\widetilde{I_j}} \arg \big( 1 + \frac{h_j}{f_j}
\big)\, \big| \lesssim R^{-C_1/2}
\end{multline*}
(we have used that $\displaystyle \Ex \Delta_{\widetilde{I_j}}
\arg T_{w_j}f = \Ex \Delta_{\widetilde{I_j}} \arg f_j$), whence
\[
\sum_{j\colon \Omega_j\ {\rm doesn't\ occur}  } \big| \delta
(T_{w_j} f, \widetilde{I_j}) - \delta (f_j, \widetilde{I_j}) \big|
\lesssim 2\pi R \cdot R^{-C_1/2} \ll 1\,.
\]
Therefore, we conclude that
\begin{multline*}
\Big| \sum_j \big(  \mathcal
\delta (T_{w_j} f, \widetilde{I_j}) - \delta (f_j, \widetilde{I_j})
\big)  \Big| \\
\le  \sum_{j\colon \Omega_j\, {\rm occurs}} \big| \delta (T_{w_j}
f, \widetilde{I_j})\big| + \sum_{j\colon \Omega_j\, {\rm occurs}}
\big| \delta (f_j, \widetilde{I_j}) \big| + 1 \\
= \sum_{j\colon \Omega_j\, {\rm occurs}} \big| \delta (f,
I_j)\big| + \sum_{j\colon \Omega_j\, {\rm occurs}} \big| \delta
(f_j, \widetilde{I_j}) \big| + 1
\end{multline*}

To estimate the size of the two sums on the right-hand side, we
introduce the (random) counter $L = \# \big\{j\colon \Omega_j\
{\rm occurs} \big\}$. Lemma~\ref{prop7.1} (applied to $\tfrac13
\Lambda$ instead of $\Lambda$) handles the case
\[
L \le \frac{b}6 \, \frac{\Lambda}{\log R} \quad \Bigl( <
\frac{b\Lambda}{3(r^2 + \log R)} \quad {\rm for} \quad R\gg 1
\Bigr).
\]
It yields that outside of some event $\Omega'$ of probability at
most $2e^{- b_1r^2 \Lambda}$, each of these two sums does not
exceed
$\frac13 \Lambda$.

Now, consider the second case when $ L > \tfrac{b}6 \,
\tfrac{\Lambda}{\log R}$. Denote by $Q$ the integer
part of $ \displaystyle \frac{b}6 \, \frac{\Lambda}{\log R}$.
Then at least $Q$ {\em independent events} $\Omega_{j_1}$, ... ,
$\Omega_{j_Q}$ must occur. By Lemma~\ref{lemma4.2} applied with
$\gamma = \widetilde{I_j}$ and  $\epsilon = R^{-C_1/2}$, we have
\[
\Pr {\Omega_j} \le 100rR^{-C_1/2}\sqrt{\tfrac12 C_1\log R} \le
R^{-C_1/3}\,,
\]
provided that $R$ is sufficiently big. Therefore,
\[
\Pr { L\ge \tfrac16\, \tfrac{b\Lambda}{\log R} } \le
\underbrace{\binom{\# \{I_j\}}{Q}}_{\le (2\pi R)^Q} \big(
R^{-C_1/3} \big)^Q  \\ < e^{-\frac14 C_1 Q \log R} \le
e^{-c_2\Lambda}\,.
\]
Thereby,
\begin{multline*}
\Pro { \big| \sum_{j} \delta(f, I_j) - \sum_j \delta(f_j,
\widetilde{I_j})  \big|  \ge \Lambda } \\ \le \Pr {\Omega'} +  \Pr
{
L\ge \tfrac16\, \tfrac{b\Lambda}{\log R} }
< 2e^{-b_1r^2 \Lambda} + e^{-c_2\Lambda} < e^{-c_3\Lambda}\,,
\end{multline*}
and we are done. \hfill $\Box$

\subsection{Proof of Theorem~\ref{thm_main}:
the upper bound in the case $\tfrac12 <\alpha<1$}

We split the circumference $R\T$ into $ N = \lfloor 2\pi R \rfloor
$ disjoint arcs $\{I_j\}$ of equal length $r$, $1\le r \le 2$. We
fix a positive $\epsilon<\tfrac{1-\alpha}4$ and suppose that
\[
\Bigl| \sum_{j=1}^N \delta (f, I_j)  \Bigr|
> 2\pi R^\alpha\,.
\]
Then we split the set $\big\{ 1, \, ... \,, N \big\}$ into $n =
\lfloor
2R^\epsilon \rfloor $ disjoint arithmetic progressions $J_1$, ...,
$J_n$. If $R$ is sufficiently big, then the cardinality of each of
these arithmetic progressions cannot be less than
\[
\frac{N}{n} - 1 \ge \frac{2\pi R -1}{2R^\epsilon} -1 >
2R^{1-\epsilon},
\]
and cannot be larger than
\[
\frac{N}{n} + 1 \le \frac{2\pi R}{2R^\epsilon -1} + 1 < 4
R^{1-\epsilon}\,.
\]
For at least one of these progressions, say for $J_l$, we have
\[
\Bigl| \sum_{j\in J_l} \delta (f, I_j) \Bigr|
> 2\pi\, \frac{R^\alpha}{n} > 2R^{\alpha-\epsilon}\,.
\]

Given a collection $\{I_j\}_{j\in J}$ with
$ 2R^{1-\epsilon} < \# J < 4 R^{1-\epsilon}  $
of $ R^\epsilon $-separated arcs of length $r$, we
show that
\[
\Pro { \Bigl| \sum_{j\in J} \delta (f, I_j)  \Bigr|
> 2R^{\alpha - \epsilon} } \le C_1 e^{-c_2 R^{2\alpha -1 -
\epsilon}}\,.
\]
Since we have  $n\ll R$ such collections $\big\{
I_j\big\}$, this will prove the upper bound in the case $\tfrac12
< \alpha < 1$.

Now, suppose that $ \bigl| \sum_{j\in J} \delta (f, I_j)  \bigr|
> 2R^{\alpha - \epsilon} $.
By Lemma~\ref{sub-thm2} applied with $\Lambda =
R^{\alpha-\epsilon}$, we see that there is a collection of
independent G.E.F. $\{f_j\}$ such that throwing away an event
of probability at most
\[ e^{-b_2 \Lambda} = e^{-b_2R^{\alpha - \epsilon }}
\stackrel{\epsilon < 1-\alpha}\ll
e^{-R^{2\alpha -1}}\,,\] we have
\[
\Bigl| \sum_{j\in J} \delta (f_j, \widetilde{I_j}) \Bigr| >
2R^{\alpha-\epsilon} - \Lambda =  R^{\alpha-\epsilon}\,.
\]
To estimate the probability of the event $ \Pr { \bigl| \sum_{j\in
J} \delta (f_j, \widetilde{I_j}) \bigr| > R^{\alpha-\epsilon} }  $,
we apply Bernstein's estimate (Lem\-ma~\ref{lemma2.3}) to the
independent identically distributed random variables
$ \psi_j =  \delta (f_j, \widetilde{I_j}) $.
By Lemma~\ref{lemma.new-old}, the tails of these random variables
decay superexponentially:
\[
\Pr { \bigl| \psi_j \bigr| \ge t } \le
\exp\Bigl(- \frac{c_3t^2}{\log t} \Bigr)
\]
for $t\gg 1$. The number of the random variables $\psi_j$ is bigger
than $ 2R^{1-\epsilon}$. Hence, the Bernstein estimate can be
applied with $ t = R^{\alpha-\epsilon} $. We see that the
probability we are interested in does not exceed \[ 2 \exp \big( -
c_4 t^2/\# J \big) < \exp \big( -c_5
R^{2\alpha-1-\epsilon}\big)\]
completing the argument. \hfill $\Box$

\section{Proof of Theorem~\ref{thm_main}: the lower bound for
$\tfrac12 <\alpha<1$}\label{sect5.3}
\numberwithin{equation}{section}

We fix $\alpha\in (\tfrac12, 1)$ and show that, for some positive
numerical constant $c_0$ and for each $R>R_0(\alpha)$, one
has
\[
\Pr { n(R) \le R^2 - c_0 R^\alpha } \ge e^{-3R^{2\alpha-1}}\,.
\]
Everywhere below, we assume that $R>2$. Let $N=\lfloor R\rfloor$.
Let $\J_-$ be a set consisting of $N$ integers between $R^2-2R$
and $R^2-R$, and let $\J_+$ be a set consisting of $N$ integers
between $R^2+R$ and $R^2+2R$. Let
$$
a_k=\begin{cases} \sqrt{1-R^{\alpha-1}},& k\in\J_+\,;
\\
\sqrt{1+R^{\alpha-1}},& k\in\J_-\,;
\\
1,& k\notin{\J_+\cup\J_-}\,.
\end{cases}
$$
Consider the Gaussian Taylor series
\[
g(z) = \sum_{k=0}^\infty \zeta_k a_k \frac{z^k}{\sqrt{k!}}\,,
\]
and denote by $n_g(R)$ the number of its zeroes in the disk $R\D$.
\begin{claim}\label{claim5.1} For $R\ge 1$, we have
$\Ex n_g(R)\le R^2 - c_1 R^\alpha$.
\end{claim}

\par\noindent{\em Proof of Claim~\ref{claim5.1}:} By
Lemma~\ref{lemma2.EK},
\[
\Ex n_g(R) = \frac12\, \frac{R \sum_{k\ge 0} a_k^2\cdot 2k \cdot
\frac{R^{2k-1}}{k!}}{\sum_{k\ge 0} a_k^2 \cdot \frac{R^{2k}}{k!}}
= \frac{\sum_{k\ge 0} a_k^2 \cdot k \cdot
\frac{R^{2k}}{k!}}{\sum_{k\ge 0} a_k^2 \cdot \frac{R^{2k}}{k!}}\,.
\]
The ratio on the right-hand side can be written as
$$
R^2 + \frac{\sum_{k\ge 0}a_k^2\cdot (k-R^2)\cdot
\frac{R^{2k}}{k!}} {\sum_{k\ge 0}a_k^2\cdot\frac{R^{2k}}{k!}}\,.
$$
Note that
$$
\sum_{k\ge 0}(k-R^2)\cdot \frac{R^{2k}}{k!}=0\,,
$$
so the numerator in the second term equals
\begin{multline*}
\sum_{k\in\J_-}R^{\alpha-1}\cdot (k-R^2)\cdot \frac{R^{2k}}{k!} +
\sum_{k\in\J_+}(-R^{\alpha-1})\cdot (k-R^2)\cdot \frac{R^{2k}}{k!}
\\
\le -R^\alpha\, \sum_{k\in\J_-\cup\J_+}\frac{R^{2k}}{k!}\,.
\end{multline*}
Since $R\ge 1$, we have  $a_k^2\le 2$, and the denominator cannot
be
bigger than $2 e^{R^2}$. Hence,
\[
\Ex n_g(R) \le R^2  - \frac12 R^\alpha e^{-R^2}
\sum_{k\in\J_-\cup\J_+}\frac{R^{2k}}{k!}\,.
\]

Now, observe that
$$
\sum_{k\in\J_-\cup\J_+}\frac{R^{2k}}{k!}\ge c e^{R^2}
$$
with some absolute $c>0$. To see this, note that the function
$k\mapsto \tfrac{R^{2k}}{k!}$ decreases for $k\in\J_+$ and
increases for $k\in\J_-$. We set $K = \lceil R^2+2R \rceil $.
Applying Stirling's formula, we get
\begin{eqnarray*}
\frac{R^{2k}}{k!} \ge \frac{R^{2K}}{K!} &\gtrsim&
\frac1{\sqrt{K}}\,\left(\frac{eR^2}{K}\right)^K   \\
&\gtrsim& \frac{e^{R^2+2R-1}}{R\left(
1+\tfrac2{R}\right)^{R^2+2R}} \gtrsim \frac{e^{R^2+2R}}{R
e^{2R+4}} \gtrsim \frac{e^{R^2}}{R}
\end{eqnarray*}
for $k\in\J_+$. A similar estimate holds for $k\in\J_-$.
Therefore,
\[
\Ex n_g(R) \le R^2  - \frac12 R^\alpha e^{-R^2}
\sum_{k\in\J_-\cup\J_+}\frac{R^{2k}}{k!} \le R^2 - \frac12
R^\alpha e^{-R^2} \cdot ce^{R^2}
\]
proving the claim. \hfill $\Box$

\begin{claim}\label{claim5.2} For $R\ge 1$, we have
$$
\Pr { n_g(R)\le R^2-\frac{c_1}2 R^{\alpha} }  \ge \frac{c_1}2
R^{-2+\alpha}\,.
$$
\end{claim}

\par\noindent{\em Proof of Claim~\ref{claim5.2}:}
We have
$$
c_1 R^{\alpha} \le \Ex (R^2-n_g(R)) \le \frac{c_1}2 R^{\alpha} +
R^2 \Pr { n_g(R)\le R^2-\frac{c_1}2 R^{\alpha}}
$$
whence
$$
\Pr {n_g(R) \le R^2-\frac{c_1}2 R^{\alpha} } \ge
R^{-2}\cdot\frac{c_1}2 R^{\alpha} = \frac{c_1}2 R^{-2+\alpha}\,.
$$
\hfill $\Box$

\begin{claim}\label{claim5.3}
Let $0\le t\le N$. Then
$$
\Pro { \sum_{k\in\J_-} |\zeta_k|^2 - \sum_{k\in\J_+}|\zeta_k|^2\ge
t } \le 2\exp\left( -\frac{t^2}{16(e+1)N} \right)\,.
$$
\end{claim}

\par\noindent{\em Proof of Claim~\ref{claim5.3}:}
Note first of all that $\Pr { |\zeta_k|^2\ge t } = e^{-t}$ and
$\Ex |\zeta_k|^2=1$, whence, for $t>0$,
\[
\Pr {|\zeta_k|^2-1 > t} < e^{-t}
\]
and
\[
\Pr {|\zeta_k|^2-1 < -t} = \max\bigl\{1 - e^{t-1}, 0 \bigr\} <
e^{1-t}\,.
\]
Thus we can apply Bernstein's lemma~\ref{lemma2.3} with $K=e+1$ to
the random variables $\pm (|\zeta_k|^2-1)$, which yields the
desired conclusion. \hfill $\Box$

\bigskip In particular,
$$
\Pro { \sum_{k\in\J_-}|\zeta_k|^2-\sum_{k\in\J_+}|\zeta_k|^2\ge
R^{1/2}\log R } \le 2\exp\left( -c_2 \log^2 R \right)\le
\frac{c_1}4 R^{-2 + \alpha}\,,
$$
provided that $R>R_0(\alpha)$.

\bigskip Now everything is ready to make the final estimate.
Let $\gamma $ be the standard Gaussian measure on the space
$\C^\infty$; i.e., the product of countably many copies of the
measures $\displaystyle \frac1{\pi} e^{-|\eta_k|^2}\, dm(\eta_k)$,
and let $\gamma_a$ be another Gaussian measure on $\C^\infty$ that
is the product of the Gaussian measures $\displaystyle \frac1{\pi
a_k^2} e^{-|\eta_k|^2/a_k^2}\, dm(\eta_k)$. Let $E\subset
\C^\infty $ be the set of coefficients $\eta_k$ such that the
Taylor series $\sum_{k\ge 0}\eta_k\frac{z^k}{\sqrt{k!}}$ converges
in $\C$ and has at most $R^2-\tfrac{c_1}2 R^{\alpha}$ zeroes in
$R\D$. Then Claim~\ref{claim5.2} can be rewritten as
\[
\gamma_a (E) \ge \frac{c_1}2 R^{-2+\alpha}\,,
\]
while the quantity $\Pr { n(R)\le R^2 - \tfrac{c_1}2 R^\alpha}$ we
are interested in equals $\gamma (E)$. Thus, it remains to compare
$\gamma (E)$ with $\gamma_a (E)$.

Let
$$
U = \Bigl\{\sum_{k\in \J_-\cup\J_+}|\eta_k|^2 \ge \sum_{k\in
\J_-\cup\J_+}\frac{|\eta_k|^2}{a_k^2} + R^{\alpha - \frac12}\log R
\Bigr\}\,.
$$
and
$$
\widetilde U = \Bigl\{\sum_{k\in \J_-\cup\J_+}a_k^2 |\eta_k|^2\ge
\sum_{k\in \J_-\cup\J_+}|\eta_k|^2 +R^{\alpha - \frac12}\log R
\Bigr\}\,.
$$
Note that
\[
\gamma_a (U) = \gamma (\widetilde{U})  = \Pro {
\sum_{k\in\J_-}|\zeta_k|^2-\sum_{k\in\J_+}|\zeta_k|^2\ge
R^{1/2}\log R } \le \frac{c_1}4 R^{-2+\alpha}\,.
\]
Hence,
$$
\gamma_a (E\setminus U) \ge \frac{c_1}4 R^{-2+\alpha}\,.
$$
But on $E\setminus U$, we can bound the density of $\gamma_a$ with
respect to $\gamma$:
$$
\frac{d\gamma_a}{d\gamma} \le e^{R^{\alpha -\frac12}\log
R}(1-R^{2\alpha-2})^{-N}  < e^{2R^{2\alpha-1}}
$$
for $R>R_0(\alpha)$. The rest is obvious:
\[
\gamma (E) \ge \gamma (E\setminus U) \ge e^{-2R^{2\alpha-1}}
\gamma_a(E\setminus U) \\
\ge  \frac{c_1}4 R^{-2+\alpha} e^{-2R^{2\alpha-1}}  \ge
e^{-3R^{2\alpha-1}},
\]
provided that $R>R_0(\alpha)$. This proves the lower bound in
Theorem~\ref{thm_main}. \hfill $\Box$

\section*{Appendix: Asymptotic almost independence.
Proof of Lemma~\ref{thm3.1}}\nonumber

\numberwithin{equation}{subsection}
\renewcommand{\thesubsection}{A-\arabic{subsection}}

% redefine the command that creates the equation no.
\setcounter{subsection}{0}  % reset counter

\subsection{Elementary inequalities}

\begin{claim}\label{claim_1a} For all positive $k$ and $t$,
\[
k \log t - t \le k\log k - k - (\sqrt{t} - \sqrt{k})^2\,.
\]
\end{claim}

\par\noindent{\em Proof:} The function $\phi (\tau) = k \log
(\tau^2) - \tau^2$ attains its maximum at $\tau  = \sqrt{k}$, and
\[
\phi''(\tau) = -\frac{2k}{\tau^2} - 2 \le -2 \qquad {\rm for\ all\
} \tau>0\,.
\]
Hence,
\[
\phi (\tau) \le \phi(\sqrt{k}) - (\tau-\sqrt{k})^2 \qquad {\rm
for\ all\ } \tau>0\,.
\]
Replacing $\tau^2$ by $t$, we get the claim. \hfill $\Box$

\begin{claim} Let $k$ be a positive integer and $u\ge k$. Then
\[
\int_u^\infty \frac{t^k e^{-t}}{k!}\, dt \le e^{-(\sqrt{u} -
\sqrt{k})^2}\,.
\]
\end{claim}

\par\noindent{\em Proof:}
\begin{multline*}
\int_u^\infty \frac{t^k e^{-t}}{k!}\, dt = \int_k^\infty
\frac{[t+(u-k)]^k e^{-t-(u-k)}}{k!}\, dt \\
= \int_k^\infty \frac{t^k e^{-t}}{k!} \left[ 1 + \frac{u-k}{t}
\right]^k e^{-(u-k)}\, dt \\
\le \left[ 1 + \frac{u-k}{k} \right]^k e^{-(u-k)}
\underbrace{\int_k^\infty \frac{t^k e^{-t}}{k!}\, dt}_{\le 1} \\
\le \exp \big\{ [k\log u - u] - [k\log k - k] \big\}
\stackrel{\rm Claim~\ref{claim_1a}}\le e^{-(\sqrt{u} - \sqrt{k})^2}
\end{multline*}
proving the claim. \hfill $\Box$

\begin{corollary}\label{claim_1c}
\[
\frac1{\pi} \int_{|z|\ge \sqrt{k}+d} \frac{|z|^{2k}}{k!}
e^{-|z|^2}\, dm_2(z) = \int_{(\sqrt{k}+d)^2}^\infty
\frac{t^k}{k!}\, e^{-t}\, dt \le e^{-d^2}\,.
\]
\end{corollary}

\begin{claim}\label{claim_1d} Let $w'$, $w''$ be points in $\C$
and let $k'$, $k''$ be non-negative integers. Then
\[
\big| \Ex \big\{ \xi_{k'}(w') \overline{\xi_{k''}(w'')} \big\}
\big| \le 2e^{-\frac{d^2}8}\,,
\]
provided that $|w'-w''|\ge \sqrt{k'} + \sqrt{k''} + d$, $d>0$.
\end{claim}

\par\noindent{\em Proof:} By Section~\ref{sect2.2}(d),
\begin{multline*}
\Bigl| \Ex \big\{ \xi_{k'}(w') \overline{\xi_{k''}(w'')} \big\}
\Bigr| = \Bigl| \Big\langle T_{-w'} \big(
\frac{z^{k'}}{\sqrt{k'!}} \big),
T_{-w''} \big( \frac{z^{k''}}{\sqrt{k''!}} \big) \Big\rangle \Bigr|
\\
= \Bigl| \frac1{\pi} \int_{\C}   \frac{(z-w')^{k'}}{\sqrt{k'!}}\,
\frac{(\overline{z-w''})^{k''}}{\sqrt{k''!}}\, e^{-z
\overline{w}'-\frac12{|w'|^2}} e^{-\overline{z} w''
-\frac12{|w''|^2}} e^{-|z|^2}\, dm_2(z) \Bigr|\,.
\end{multline*}
Therefore,
\begin{multline*}
\Big| \Ex \big\{ \xi_{k'}(w') \overline{\xi_{k''}(w'')} \big\}
\Big| \\
\le \frac1{\pi} \int_{\C} \frac{|z-w'|^{k'}}{\sqrt{k'!}}\,
e^{-\frac12 |z-w'|^2} \, \frac{|z-w''|^{k''}}{\sqrt{k''!}}\,
e^{-\frac12 |z-w''|^2} \, dm_2(z) \\
\le \frac1{\pi}\int_{\big\{ |z-w'|\ge \sqrt{k'}+\frac{d}2\big\}}
\, +\,  \frac1{\pi}\int_{\big\{ |z-w''|\ge
\sqrt{k''}+\frac{d}2\big\}}  = I' + I''\,.
\end{multline*}
By the Cauchy-Schwarz inequality,
\begin{multline*}
I' \le \Big\{ \frac1{\pi} \int_{\big\{ |z-w'|\ge
\sqrt{k'}+\frac{d}2\big\}} \frac{|z-w'|^{2k'}}{k'!}\,
e^{- |z-w'|^2} \, dm_2(z)  \Big\}^{1/2} \\
\times \Big\{ \frac1{\pi} \int_{\C} \frac{|z-w''|^{2k''}}{k''!}\,
e^{- |z-w''|^2} \, dm_2(z) \Big\}^{1/2} \\
\stackrel{\rm Claim~\ref{claim_1c}}\le \, e^{-\frac{d^2}8} \cdot 1
=
e^{-\frac{d^2}8}\,.
\end{multline*}
Similarly, $\displaystyle I'' \le e^{-\frac{d^2}8}$. Hence,
$\displaystyle I'+I''\le 2 e^{-\frac{d^2}8}$, and we are done.
\hfill $\Box$

\begin{claim}\label{claim_1e} Assume that the disks $D(w_j,
R_j+8\sigma_j)$ are pairwise disjoint and $R_j\ge 1$, $\sigma_j
\ge \max \big( 1, \sqrt{\log R_j}\big)$. Let $D_{ij} =
|w_i-w_j|-R_i-R_j$ be the distance between the disks $D(w_i, R_i)$
and $D(w_j, R_j)$. Then, for each $i$,
\[
2\sum_{j\colon\, j\ne i} (1+R_j^2) e^{-\frac18 D_{ij}^2} \le
e^{-2\sigma_i^2}\,.
\]
\end{claim}

\par\noindent{\em Proof:} Indeed, since $D_{ij}\ge 8\sigma_j$, we
have
\[
\frac1{16} D_{ij}^2 \ge 4 \sigma_j^2 \ge 2\sigma_j^2 +2 \ge \log
(e^2 R_j^2)  \ge \log (4R_j^2) \ge \log\big[ 2(1+R_j^2)\big]\,.
\]
Thus, it suffices to estimate the sum $\displaystyle \sum_{j\colon
j\ne i} e^{-\frac1{16} D_{ij}^2}$. For each $j\ne i$, consider the
disk $\mathcal D_j \subset D(w_j, R_j+8\sigma_j)$ of radius $4$
closest to $w_i$.
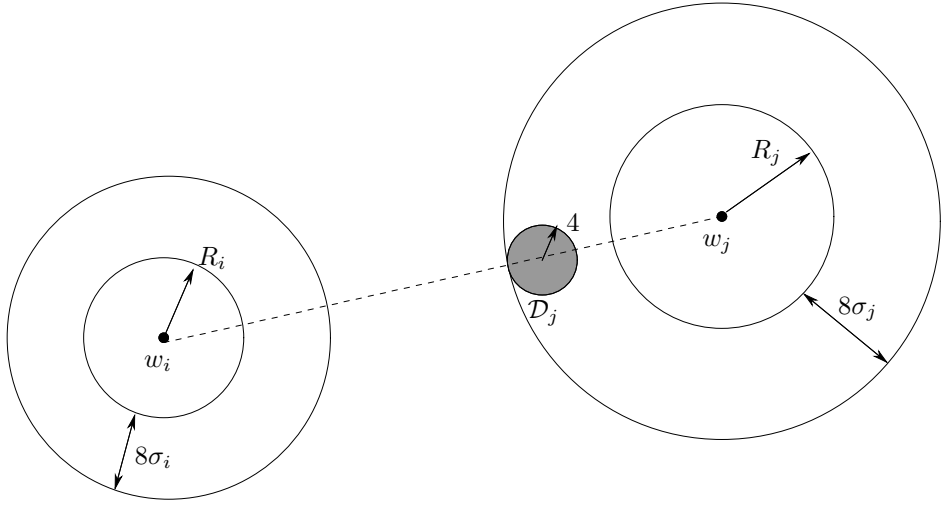
\begin{figure}[h]
\begin{center}
\input{disks.eepic}
\end{center}
\caption{\tt The disks $D(w_i, R_i)$, $D(w_j, R_j)$ and  $\mathcal
D_j$}
\end{figure}
For each $z\in\mathcal D_j$, we have $|z-w_i|\le D_{ij}+R_i$.
Also, the disks $\mathcal D_j$ are disjoint and $\displaystyle
\bigcup_j \mathcal D_j \subset \C \setminus D(w_i,
R_i+8\sigma_i)$. Hence,
\begin{multline*}
\sum_{j\colon j\ne i} e^{-\frac1{16} D_{ij}^2} \le \frac1{16\pi}
\int_{\{ |z-w_i|\ge R_i+8\sigma_i \}} e^{-\frac1{16}
(|z-w_i|-R_i)^2 }\, dm_2 (z)\\
= \frac1{16\pi} \int_{\{|z|\ge R_i+8\sigma_i\}} e^{-\frac1{16}
(|z|-R_i)^2}\, dm_2(z)  = \frac18 \int_{8\sigma_i}^\infty
(R_i+t)e^{-\frac1{16}t^2}\, dt
\\
\le (1 + \tfrac18 R_i) \int_{8\sigma_i} \frac{t}8 e^{-\frac1{16}
t^2}\, dt =
(1 + \tfrac18 R_i) e^{-4\sigma_i^2} \\
\le \frac{1 + \tfrac18 R_i}{e^{\sigma_i^2+1}} e^{-2\sigma_i^2} \le
\frac{8 + R_i}{8eR_i} e^{-2\sigma_i^2} \le \frac9{8e}
e^{-2\sigma_i^2} < e^{-2\sigma_i^2}
\end{multline*}
proving the claim. \hfill $\Box$

\subsection{Almost orthogonal standard Gaussian random variables
are almost independent}

\begin{claim}\label{claim_2}
Let $\xi_j$ be standard complex Gaussian random variables such
that their covariance matrix $\Gamma_{ij} = \Ex\big\{ \xi_i
\overline{\xi_j} \big\}$ satisfies
\[
\sum_{j\colon j\ne i} |\Gamma_{ij}| \le \delta_i \le \frac13\,.
\]
Then $\xi_j = \zeta_j + s_j \eta_j$ where $\zeta_j$ are
independent standard complex Gaussian random variables, $\eta_j$
are standard complex Gaussian random variables, and $s_j\in [0,
\delta_j]$.
\end{claim}

\par\noindent{\em Proof:} Let $\Gamma = I-\Delta$ where $I$ is the
identity matrix. Put
\[
\zeta_i = \sum_j (\Gamma^{-1/2})_{ij} \xi_j\,.
\]
Then $\zeta_i$ are independent standard complex Gaussian random
variables. We set $\widetilde{\Delta} = I - \Gamma^{-1/2}$ and
$\displaystyle s_i \eta_i = \sum_j \widetilde{\Delta}_{ij} \xi_j$,
and estimate the sum $\displaystyle \sum_j
|\widetilde{\Delta}_{ij}|$.

We have
\[
\Gamma^{-1/2} = I + \frac12 \Delta + \sum_{k\ge 2} \alpha_k
\Delta^k
\]
with $|\alpha_k|\le 1 $ for all $k\ge 2$. Then
\[
|\widetilde{\Delta}_{ij}| \le   \frac12 |\Delta_{ij}| + \sum_{k\ge
2} |(\Delta^k)_{ij}|,
\]
whence
\[
\sum_j |\widetilde{\Delta}_{ij}| \le \frac12 \sum_{j}
|\Delta_{ij}| + \sum_{k\ge 2} \sum_j |(\Delta^k)_{ij}| \le
\frac{\delta_i}2 + \sum_{k\ge 2} \sum_j |(\Delta^k)_{ij}| \,.
\]

To estimate the sum on the right-hand side, we note that for any
two square matrices $A$ and $B$ of the same size, we have
\begin{multline*}
\sum_j |(AB)_{ij}| \le \sum_{j, \ell} |A_{i \ell}|\, |B_{\ell j}|
\\
= \sum_{\ell} \big[ |A_{i\ell}| \cdot \sum_j |B_{\ell j}| \big]
\le \big( \sum_j |A_{ij}|\big) \cdot \sup_{\ell} \sum_j |B_{\ell
j}|\,.
\end{multline*}
Applying this observation to the matrices $\Delta^k = \Delta \cdot
\Delta^{k-1} $ (with $k\ge 1$), we conclude by induction that
\[
\sum_{j} |(\Delta^k)_{ij}| \le \big( \sum_j |\Delta_{ij}| \big)
3^{-(k-1)} \le \frac{\delta_i}{3^{k-1}}\,.
\]
Thus
\[
\sum_j |\widetilde{\Delta}_{ij}| \le \frac{\delta_i}2 + \sum_{k\ge
2} \frac{\delta_i}{3^{k-1}} = \delta_i\,,
\]
and we are done. \hfill $\Box$

\subsection{Proof of the lemma}

We fix two big constants $A \gg a \gg 1$.
Let $R_j = r_j + a\rho_j$, $\sigma_j = \frac{A-a}8 \rho_j$.
Clearly, $R_j\ge 1$, $\sigma_j \ge 1$. Also,
\begin{multline*}
\sigma_j = 2\rho_j + \left( \frac{A-a}8 - 2 \right) \rho_j \\
\ge 2\sqrt{\log r_j} + \frac{A-a-16}{8a} \log (1+a\rho_j) \\
\ge 2 \sqrt{\log r_j} + 2\sqrt{\log(1+a\rho_j)} \\
\ge 2\sqrt{\log r_j (1+a\rho_j)} \ge 2\sqrt{\log R_j}\,,
\end{multline*}
provided that $a\ge 2$ and $A \ge 17a + 16$.

We consider now the family of standard Gaussian random variables
$\zeta_k(w_j)$, $k\le R_j^2$. Applying  to this family
Claim~\ref{claim_1d}, we get
\[
\big| \Ex \big\{ \zeta_k (w_i) \overline{\zeta_{\ell}(w_j)} \big\}
\big| \le 2e^{-\frac18 D_{ij}^2}
\]
where, as before, $D_{ij} = |w_i-w_j| - R_i - R_j$ is the distance
between the disks $D(w_i, R_i)$ and $D(w_j, R_j)$. Now,
Claim~\ref{claim_1e} implies that the sum of absolute values of
the covariances of $\zeta_k(w_i)$ with all other $\zeta_l (w_j)$
in our family does not exceed $\displaystyle e^{-2\sigma_i^2} \le
e^{-2} < \frac13$. Claim~\ref{claim_2} then allows us to write
\[
\zeta_k (w_i) = \zeta_{ik} + s_{ik} \eta_{ik}\,, \qquad k\le R_i^2
\]
where $\zeta_{ik}$ are independent standard Gaussian complex
random variables, $\eta_{ik}$ are standard Gaussian complex random
variables, and $s_{ik}\in [0, e^{-2\sigma_i^2}]$.

Next, we choose $\zeta_{ik}$, $k>R_i^2$, in such a way that the
whole family $\zeta_{ik}$ of standard Gaussian complex random
variables is independent and put
\begin{eqnarray*}
f_i &=& \sum_k \zeta_{ik} \frac{z^k}{\sqrt{k!}}\,, \\
h_i &=& \sum_{k\le R_i^2} s_{ik} \eta_{ik} \frac{z^k}{\sqrt{k!}} +
\sum_{k>R_i^2} \big[ \zeta_k(w_i) - \zeta_{ik} \big]
\frac{z^k}{\sqrt{k!}}\,.
\end{eqnarray*}
By construction, $T_{w_i} f = f_i + h_i$.

To estimate the probability
\[
\Pr { \max_{z\in r_i\D} |h_i(z)|e^{-\frac12 |z|^2} >
e^{-\rho_j^2} }\,,
\]
it suffices to estimate the expression
\[
\sum_{k\le R_i^2} s_{ik} \max_{z\in r_i\D} \frac{|z|^k}{\sqrt{k!}}
e^{-\frac12 |z|^2} + 2\sum_{k>R_i^2} \max_{z\in r_i\D}
\frac{|z|^k}{\sqrt{k!}} e^{-\frac12 |z|^2}\,.
\]
If this expression is less than $e^{-2\rho_j^2}$, then by
Lemma~\ref{lemma2.1}, we get what Lemma~\ref{thm3.1} asserts:
\[
\Pr {\max_{z\in r_j \D}  |h_j(z)| e^{-|z|^2/2} \ge e^{-\rho_j^2}
} \le 2 \exp \big( - \frac12 e^{2\rho_j^2} \big)\,.
\]

For every $k\ge 1$, we have $\displaystyle \frac{|z|^k}{\sqrt{k!}}
e^{-|z|^2/2} \le 1 $ and thereby,
\begin{multline*}
\sum_{k\le R_i^2} s_{ik} \max_{z\in r_i\D} \frac{|z|^k}{\sqrt{k!}}
e^{-|z|^2/2} \le \sum_{k\le R_i^2} s_{ik} \\
\le (1+R_i^2) e^{-2\sigma_i^2} \le
\frac{1+R_i^2}{e^{\sigma_i^2}}\, e^{-\sigma_i^2} \\
\le \frac{1+R_i^2}{R_i^4} e^{-\sigma_i^2} \le 2
e^{-\frac{(A-a)^2}{64} \rho_i^2} \le \frac12 e^{-2\rho_i^2}\,,
\end{multline*}
provided that $ A>a+16 $.

For $k>R_i^2$, Claim~\ref{claim_1a} implies that
\[
\frac{|z|^{2k}}{k!} e^{-|z|^2} \le \frac{k^k}{k!} e^{-k}\,
e^{-(\sqrt{k} - |z|)^2} \le e^{-(\sqrt{k} - r_i)^2}
\]
for all $z\in r_i\D$. Hence,
\[
\max_{z\in r_i\D} \frac{|z|^k}{\sqrt{k!}} e^{-|z|^2/2} \le
e^{-\frac12 (\sqrt{k}-r_i)^2}\,, \qquad k>R_i^2\,,
\]
and it suffices to show that
\[
2\sum_{k>R_i^2} e^{-\frac12 (\sqrt{k}-r_i)^2} \le \frac12
e^{-2\rho_i^2}\,.
\]
Now,
\[
\sum_{k>R_i^2} = \sum_{R_i^2 < k \le 4r_i^2} + \sum_{k>
\max(R_i^2, 4r_i^2)}
\]
with the usual convention that the sum taken over the empty set
equals zero. The first sum does not exceed
\[
(1+4r_i^2) e^{-\frac12 a^2\rho_i^2} \le
\frac{5r_i^2}{e^{4+2\rho_i^2}} e^{(-\frac12 a^2 -6)\rho_i^2} \,
\stackrel{\rho_i^2 \ge \log r_i}\le \, \frac5{e^4}
e^{-2\rho_i^2} < \frac18 e^{-2\rho_i^2}\,,
\]
provided that $a\ge 4$. At last, the remaining sum
does not exceed
\begin{eqnarray*}
\sum_{k\ge a\rho_i^2} e^{-\frac18 k} &\le& \frac1{1-e^{-1/8}}
e^{-\frac18 a^2\rho_i^2}  \\
&\le& 9 e^{-\frac18 a^2\rho_i^2} \le \frac9{e^6} e^{-(\frac18
a^2 - 6 ) \rho_i^2} \le \frac18 e^{-2\rho_i^2}\,,
\end{eqnarray*}
provided that $a\ge 8$. This finishes off the
proof of Lemma~\ref{thm3.1}. \hfill $\Box$

\end{document}

%% file: disks.eepic
\setlength{\unitlength}{0.00056868in}
\begingroup\makeatletter\ifx\SetFigFont\undefined%
\gdef\SetFigFont#1#2#3#4#5{%
  \reset@font\fontsize{#1}{#2pt}%
  \fontfamily{#3}\fontseries{#4}\fontshape{#5}%
  \selectfont}%
\fi\endgroup%
{\renewcommand{\dashlinestretch}{30}
\begin{picture}(8824,4635)(0,-10)
\put(1462,1507){\blacken\ellipse{90}{90}}
\put(1462,1507){\ellipse{90}{90}}
\put(1462,1507){\ellipse{1484}{1484}}
\put(6637,2632){\ellipse{2078}{2078}}
\put(6637,2632){\blacken\ellipse{90}{90}}
\put(6637,2632){\ellipse{90}{90}}
\put(1507,1507){\ellipse{2998}{2998}}
\put(6637,2587){\ellipse{4050}{4050}}
\texture{55888888 88555555 5522a222 a2555555 55888888 88555555 552a2a2a 2a555555 
	55888888 88555555 55a222a2 22555555 55888888 88555555 552a2a2a 2a555555 
	55888888 88555555 5522a222 a2555555 55888888 88555555 552a2a2a 2a555555 
	55888888 88555555 55a222a2 22555555 55888888 88555555 552a2a2a 2a555555 }
\put(4972,2227){\shade\ellipse{648}{648}}
\put(4972,2227){\ellipse{648}{648}}
\dashline{60.000}(1462,1462)(6637,2632)
\path(1462,1507)(1732,2137)
\path(1462,1507)(1732,2137)
\blacken\path(1712.304,2014.885)(1732.000,2137.000)(1657.155,2038.520)(1698.911,2059.792)(1712.304,2014.885)
\blacken\path(1190.068,663.322)(1192.000,787.000)(1132.093,678.782)(1170.356,705.836)(1190.068,663.322)
\path(1192,787)(1012,112)
\path(1192,787)(1012,112)
\blacken\path(1013.932,235.678)(1012.000,112.000)(1071.907,220.218)(1033.644,193.164)(1013.932,235.678)
\blacken\path(7513.703,1858.873)(7402.000,1912.000)(7475.560,1812.557)(7466.842,1858.601)(7513.703,1858.873)
\path(7402,1912)(8167,1282)
\path(7402,1912)(8167,1282)
\blacken\path(8055.297,1335.127)(8167.000,1282.000)(8093.440,1381.443)(8102.158,1335.399)(8055.297,1335.127)
\path(6682,2677)(7447,3217)
\path(6682,2677)(7447,3217)
\blacken\path(7366.264,3123.289)(7447.000,3217.000)(7331.663,3172.307)(7378.375,3168.559)(7366.264,3123.289)
\path(4972,2227)(5107,2542)
\path(4972,2227)(5107,2542)
\blacken\path(5087.304,2419.885)(5107.000,2542.000)(5032.155,2443.520)(5073.911,2464.792)(5087.304,2419.885)
\put(1192,337){\makebox(0,0)[lb]{\smash{{\SetFigFont{10}{12.0}{\rmdefault}{\mddefault}
{\updefault}$8\sigma_i$}}}}
\put(1282,1237){\makebox(0,0)[lb]{\smash{{\SetFigFont{10}{12.0}{\rmdefault}{\mddefault}
{\updefault}$w_i$}}}}
\put(1777,2182){\makebox(0,0)[lb]{\smash{{\SetFigFont{10}{12.0}{\rmdefault}{\mddefault}
{\updefault}$R_i$}}}}
\put(5197,2497){\makebox(0,0)[lb]{\smash{{\SetFigFont{10}{12.0}{\rmdefault}{\mddefault}
{\updefault}$4$}}}}
\put(4837,1687){\makebox(0,0)[lb]{\smash{{\SetFigFont{10}{12.0}{\rmdefault}{\mddefault}
{\updefault}$\mathcal D_j$}}}}
\put(6457,2362){\makebox(0,0)[lb]{\smash{{\SetFigFont{10}{12.0}{\rmdefault}{\mddefault}
{\updefault}$w_j$}}}}
\put(6907,3172){\makebox(0,0)[lb]{\smash{{\SetFigFont{10}{12.0}{\rmdefault}{\mddefault}
{\updefault}$R_j$}}}}
\put(7717,1732){\makebox(0,0)[lb]{\smash{{\SetFigFont{10}{12.0}{\rmdefault}{\mddefault}
{\updefault}$8\sigma_j$}}}}
\end{picture}
}